\documentclass{article}
\usepackage{amsfonts,color}
\usepackage{amsmath}
 \usepackage{epsfig}
\usepackage{lscape}
\usepackage{amssymb}
\usepackage{bbm}

\renewcommand{\epsilon}{\varepsilon}

\DeclareMathOperator{\Div}{Div}

\DeclareMathOperator{\polar}{polar}
\newcommand{\id}{{\boldsymbol{\mathbbm{1}}}}

\allowdisplaybreaks[1]

\makeindex

\def\barr{\begin{array}}

\def\earr{\end{array}}
\def\bec#1{\begin{equation}\label{#1}}
\def\becn{\begin{equation*}}
\def\endec{\end{equation}}
\def\endecn{\end{equation*}}
\def\dd{\displaystyle}

 \usepackage{footnote}

\makeatletter
\let\@fnsymbol\@arabic
\makeatother

\begin{document}
\title{Derivation of a refined 6-parameter shell model: Descent from the three-dimensional Cosserat elasticity using a method of classical shell theory}

\author{ Mircea B\^irsan\thanks{Mircea B\^irsan, \ \  Lehrstuhl f\"{u}r Nichtlineare Analysis und Modellierung, Fakult\"{a}t f\"{u}r Mathematik,
Universit\"{a}t Duisburg-Essen, Thea-Leymann Str. 9, 45127 Essen, Germany; and  Department of Mathematics, Alexandru Ioan Cuza University of Ia\c si,  Blvd.
Carol I, no. 11, 700506 Ia\c si,
Romania;  email: mircea.birsan@uni-due.de}  
}

\maketitle

\begin{center}
\thanks{\textit{Dedicated  to Professor Sanda Cleja-\c Tigoiu on the occasion of her 70th birthday}}
\end{center}

\begin{abstract}
Starting from the three-dimensional Cosserat elasticity, we derive a two-dimensional model for isotropic elastic shells. For the dimensional reduction, we employ a derivation method similar to that used in classical shell theory, as presented systematically by Steigmann in [J. Elast. \textbf{111}: 91-107, 2013]. As a result, we obtain a geometrically nonlinear Cosserat shell model with a specific form of the strain-energy density, which has a simple expression with coefficients depending on the initial curvature tensor and on three-dimensional material constants. The explicit forms of the stress-strain relations and the local equilibrium equations are also recorded. Finally, we compare our results with other 6-parameter shell models and discuss the relation to the classical Koiter shell model.
\end{abstract}
\bigskip


\bigskip

\section{Introduction}\label{Intro}

Elastic shell theory is an important branch of the mechanics of deformable bodies, in view of its applications in engineering. It is also a current domain of active research, because scientists are looking for new shell models, with better properties. This task is not easy, since the shell model should be simple enough to be manageable in practical engineering problems, but on the other side it should be complex enough to account for relevant curvature and three-dimensional effects.    

The classical shell theory, also called the first order approximation theory, presents relatively simple shell models (e.g., the well-known Koiter shell model), but it is not applicable to all shell problems. The classical approach can be employed only if the Kirchhoff-Love hypotheses are satisfied; moreover, one can observe the effect of accuracy loss in classical shell theory for certain problems (see, e.g. \cite{Berdic-Mis92}). Therefore, more refined shell theories are needed.

One of the most general theories of shells, which has been much developed in the last decades, is the so-called 6-parameter shell theory. This approach has been initially proposed by Reissner \cite{Reissner74}. The theory of 6-parameter shells, presented in the books \cite{Libai98,Pietraszkiewicz-book04}, involves two independent kinematic fields:
the translation vector (3 degrees of freedom) and the rotation tensor (3 additional degrees of freedom).
Some of the achievements of this general shell theory have been presented in \cite{Pietraszkiewicz10,Eremeyev11,Pietraszkiewicz11}. We mention that the kinematical structure of 6-parameter shells is identical to the kinematical structure of Cosserat shells, which are regarded as deformable surfaces with a triad of rigid directors describing the orientation of material points. Thus, the rotation tensor in the 6-parameter model accounts for the orientation change of the triad of directors. General results concerning the existence of minimizers in the 6-parameter shell theory have been presented in \cite{Birsan-Neff-MMS-2014}.

In order to be useful in practice, the shell model should present a concrete (specific) form of the constitutive relations and strain-energy density. The specific form should satisfy these two requirements: the coefficients of the strain-energy density should be determined in terms of the three-dimensional material constants and they should depend on the (initial) curvature tensor $ \boldsymbol b $ of the reference configuration. In the literature of 6-parameter shells, we were not able to find a satisfactory strain-energy density for isotropic shells: the available specific forms are either too simple (in the sense that the coefficients are constant, i.e. independent of the initial curvature $ \boldsymbol b $), or they are general functions of the strain measures, which coefficients are not identified in terms of three-dimensional material constants.

Our present work aims to fill this gap and establishes a specific form for the strain-energy density of isotropic 6-parameter (Cosserat) elastic shells, together with explicit stress-strain relations, which fulfill the above requirements. In this model, we retain the terms up to the order $ O(h^3) $ with respect to the shell thickness $ h $ and derive a relatively simple expression of the strain-energy density, which can be used in applications. To obtain the two-dimensional strain-energy density (i.e., written as a function of $ (x_1,x_2) $, the surface curvilinear coordinates), we descend from a Cosserat three-dimensional elastic model and apply the derivation method from the classical theory of shells, which was systematically presented by Steigmann in \cite{Steigmann08,Steigmann12,Steigmann13}. Thus, in Section \ref{Sect2} we introduce the three-dimensional Cosserat continuum in curvilinear coordinates, with the appropriate strain measures \eqref{f1}, \eqref{f2},  equilibrium equations \eqref{f3} and constitutive relations \eqref{f4}-\eqref{f7}.
In Section \ref{Sect3},  we describe briefly the geometry of surfaces and the kinematics of 6-parameter shells, and define the shell strain tensor and bending-curvature tensor \eqref{f31}.

In the main Section \ref{Sect4},  we derive the two-dimensional shell model by performing the integration over the thickness and using the aforementioned derivation method \cite{Steigmann13}, inspired by the classical shell theory. Here, we adopt some assumptions which are common in the shell approaches (such as, for instance, the stress vectors on the major faces of the shells are of order $ O(h^3) $)
and are able to neglect some higher-order terms to obtain a simplified form of the strain-energy density \eqref{f61}.
For the sake of completeness, we also present the equilibrium equations for 6-parameter (Cosserat) shells \eqref{f82}, which we deduce from the condition that the equilibrium state is a stationary point of the energy functional.

Section \ref{Sect5} is devoted to further remarks and comments on the derived Cosserat shell model. We introduce the fourth-order tensor of elastic moduli for shells \eqref{f89}, \eqref{f93} and present the explicit form of the stress-strain relations \eqref{f100}. 
In order to compare our results with other 6-parameter shell models, we write the strain-energy density in an alternative useful form \eqref{f107}. We pay special attention to the comparison with the Cosserat shell model of order $ O(h^5) $ which has been presented recently in \cite{Birsan-Neff-MMS-2019}. Although the derivation methods are different, we obtain the same form of the strain-energy density, except for the coefficients of the transverse shear energy, which are unequal. The value of the transverse shear coefficient derived in the present work is confirmed by the results obtained previously through $ \Gamma $-convergence in  \cite{Neff_Hong_Reissner08} for the case of plates.

Finally, we discuss in Subsection \ref{Sect5.3} the relation between our 6-parameter shell model and the classical Koiter model. We show that, if we adopt appropriate restrictions (the material is a Cauchy continuum and the Kirchhoff-Love hypotheses are satisfied), we are able to reduce the form of our strain-energy density to obtain the classical Koiter energy, see \eqref{f125}.

\subsection*{Notations}\label{Not}

Let us present next some useful notations which will be used throughout this paper. The Latin indices $ i,j,k,... $ range over the set $ \{1,2,3\} $, while the Greek indices $ \alpha,\beta,\gamma,... $ range over the set $ \{1,2\} $. The Einstein summation convention over repeated indices is used.  A subscript comma preceding an index $ i $ (or $ \alpha $) designates partial differentiation with respect to the variable $ x_i $ (oder $ x_\alpha\, $, respectively), e.g. $ f,_i = \dfrac{\partial f}{\partial x_i}\; $.
We denote by $ \,\delta_i^j\, $ the Kronecker symbol, i.e. $ \,\delta_i^j=1 $ for $ i=j $, while $\,\delta_i^j=0 $ for $ i\neq j $.

We employ the direct tensor notation. Thus, $ \otimes $ designates the dyadic product, $ \id_3 = \boldsymbol g_i\otimes\, \boldsymbol g^i\, $ is the unit second order tensor in the 3-space, and $ \mathrm{axl}(\boldsymbol W) $ stands for the axial vector of any skew-symmetric tensor $ \boldsymbol W $.

Let $  \mathrm{tr}(\boldsymbol X)  $ denote the trace of any second order tensor $ \boldsymbol X $. The symmetric part, skew-symmetric part, and deviatoric part of $ \boldsymbol X $ are defined by
\[  
\mathrm{sym}\,  \boldsymbol X = \dfrac12\big( \boldsymbol X + \boldsymbol X^T\big),\qquad 
\mathrm{skew}\,  \boldsymbol X = \dfrac12\big( \boldsymbol X - \boldsymbol X^T\big),\qquad 
\mathrm{dev}_3  \boldsymbol X = \boldsymbol X - \dfrac13\, \big(\mathrm{tr}\,\boldsymbol X\big)\,\id_3\,.
\]
The scalar product between any second order tensors $ \,\boldsymbol A = A^{ij}\boldsymbol g_i\otimes \boldsymbol g_j = A_{ij}\,\boldsymbol g^i\otimes \boldsymbol g^j\,$ and $\, \boldsymbol B = B^{kl}\boldsymbol g_k\otimes \boldsymbol g_l = B_{kl}\,\boldsymbol g^k\otimes \boldsymbol g^l\,$
is denoted by
\[  
 \boldsymbol A : \boldsymbol B = \mathrm{tr}\big(\boldsymbol A^T \boldsymbol B \big) = A^{ij} B_{ij} = A_{kl} B^{kl}\,.
\]
If $ \, \underline{\boldsymbol C} = C^{ijkl}\boldsymbol g_i\otimes
\boldsymbol g_j\otimes \boldsymbol g_k\otimes \boldsymbol g_l \, $
is a fourth-order tensor, then we use the corresponding notations
\[  
\underline{\boldsymbol C} : \boldsymbol B =   C^{ijkl}B_{kl}\,\boldsymbol g_i\otimes
\boldsymbol g_j\,,\qquad 
 \boldsymbol A : \underline{\boldsymbol C} =  C^{ijkl}  A_{ij}\,
  \boldsymbol g_k\otimes \boldsymbol g_l \,,\qquad 
   \boldsymbol A : \underline{\boldsymbol C}   : \boldsymbol B =   C^{ijkl} A_{ij}\,B_{kl}\,.
\]
For any vector $ \boldsymbol v = v^i \boldsymbol g_i = v_i\, \boldsymbol g^i$ we write as usual
\[  
\boldsymbol A \boldsymbol v = A^{ij}v_j\, \boldsymbol g_i = A_{ij}v^j\, \boldsymbol g^i
\qquad\mathrm{and}\qquad
\boldsymbol v \boldsymbol A = \boldsymbol A^T \boldsymbol v 
= A^{ij}v_i\, \boldsymbol g_j = A_{ij}v^i\, \boldsymbol g^j\,.
\]

\section{Three-dimensional Cosserat elastic continua}\label{Sect2}

Let us consider a three-dimensional Cosserat body which occupies the domain 
$\Omega_\xi\subset\mathbb{R}^3$ in its reference configuration. 
The deformation is characterized by the vectorial map 
$ \boldsymbol \varphi_\xi:\Omega_\xi \rightarrow\Omega_c $ (here is $ \Omega_c\subset\mathbb{R}^3 $ the deformed configuration) and the microrotation tensor $ \boldsymbol{R}_\xi:\Omega_\xi \rightarrow \mathrm{SO}(3) $ (the special orthogonal group).

On the reference configuration $ \Omega_\xi $ we consider a system of curvilinear coordinates $(x_1,x_2,x_3)$, which are induced by the parametric representation $ \boldsymbol\Theta:\Omega_h \rightarrow\Omega_\xi\, $ with $ (x_1,x_2,x_3)\in \Omega_h\, $. Using the common notations, we introduce the covariant base vectors $\boldsymbol g_i :  =\dfrac{\partial\boldsymbol\Theta }{\partial x_i}\,= \boldsymbol\Theta,_i$ and the contravariant base vectors $\boldsymbol g^i$ with $ \boldsymbol g^j\cdot \boldsymbol g_i=\delta^j_i\, $.

Let 
$$ \boldsymbol \varphi :\Omega_h \rightarrow\Omega_c\,, \qquad \boldsymbol \varphi(x_1,x_2,x_3): = \boldsymbol \varphi_\xi\big( \boldsymbol\Theta(x_1,x_2,x_3)\big) ,$$ 
be the \textit{deformation function} and 
$$ \boldsymbol F_\xi  = \boldsymbol \varphi,_i\otimes\, \boldsymbol g^i\, $$ 
the \textit{deformation gradient}. We refer the domain $ \Omega_h $ to the orthonormal vector basis $ \{\boldsymbol e_1, \boldsymbol e_2 , \boldsymbol e_3 \} $, such that $ (x_1,x_2,x_3)= x_i\boldsymbol e_i\, $ and $   \nabla_x \boldsymbol\Theta=\boldsymbol\Theta,_i\,\otimes \boldsymbol e_i=\boldsymbol g_i\otimes \boldsymbol e_i$\,. The microrotation tensor can be represented as 
\[ \boldsymbol{R}_\xi = \boldsymbol d_i \otimes \boldsymbol d_i^0\,,\]
where $ \{ \boldsymbol d_1^0\,,  \boldsymbol d_2^0\,, \boldsymbol d_3^0\, \} $ is the orthonormal triad of directors in the reference configuration $ \Omega_\xi $ and $ \{ \boldsymbol d_1\,,  \boldsymbol d_2\,, \boldsymbol d_3\, \} $ is the orthonormal triad of directors in the deformed configuration $ \Omega_c\, $. We denote by $\boldsymbol Q_e$ the \emph{elastic microrotation} given by
$$
\boldsymbol Q_e :\Omega_h \rightarrow \mathrm{SO}(3 ),\qquad \boldsymbol Q_e (x_1,x_2,x_3):= \boldsymbol R_\xi\big(\boldsymbol\Theta(x_1,x_2,x_3)\big).
$$
We choose the initial microrotation tensor $ \boldsymbol Q_0 $ such that
\begin{equation}\label{f0,5}
{\boldsymbol Q}_0=\polar{(\nabla_x \boldsymbol \Theta)}\in \rm{SO}(3 )\qquad\mbox{and}\qquad \boldsymbol Q_0=\boldsymbol d_i^0\otimes\boldsymbol e_i\,.
\end{equation}
Let 
\begin{equation}\label{f1}
 \overline{\boldsymbol{E}}:=\boldsymbol Q_e^T \boldsymbol F_\xi -\id_3
\end{equation}
denote the (non-symmetric) \emph{strain tensor} for nonlinear micropolar media and 
\begin{equation}\label{f2}
\boldsymbol \Gamma := \mathrm{axl}\big(\boldsymbol Q_e^T\boldsymbol Q_{e,i}\big)\otimes \boldsymbol g^i
\end{equation}
be the \emph{wryness tensor} (see e.g., \cite{Neff_curl08,Pietraszkiewicz09,Birsan-Neff-L58-2017}), which is a strain measure for curvature (orientation change). 

The local equations of equilibrium can be written in the form 
\begin{equation}\label{f3}
\Div \boldsymbol T +\boldsymbol f=\boldsymbol 0,\qquad \Div \overline{\boldsymbol M} - \mathrm{axl}\big(\boldsymbol F_\xi \boldsymbol T^T - \boldsymbol T^T \boldsymbol F_\xi \big) + \boldsymbol c = \boldsymbol 0,
\end{equation}
where $ \boldsymbol T $ and $ \overline{\boldsymbol M} $ are the stress tensor and the couple stress tensor (of the first Piola-Kirchhoff type), $ \boldsymbol f $ and $ \boldsymbol c $ are the external body force and couple vectors. To the balance equations (\ref{f3}) one can adjoin boundary conditions.

Under hyperelasticity assumptions, the stress tensors $ \boldsymbol T $ and $ \overline{\boldsymbol M} $  are expressed by the constitutive equations 
\begin{equation}\label{f4}
\boldsymbol Q_e^T \boldsymbol T = \dfrac{\partial W}{\partial\overline{\boldsymbol{E}}}\;
,\qquad 
\boldsymbol Q_e^T \overline{\boldsymbol M}  = \dfrac{\partial W}{\partial\boldsymbol \Gamma}\;,
\end{equation}
where $ W=W(\overline{\boldsymbol{E}}, \boldsymbol \Gamma ) $ is the elastically stored energy density. Using the Cosserat model for isotropic materials presented in \cite{Birsan-Neff-Ost_L56-2015,Birsan-Neff-MMS-2019}, we assume the following representation for the energy density
\begin{align}
W(\overline{\boldsymbol{E}}, \boldsymbol \Gamma)=\; & W_{\mathrm{mp}}(\overline{\boldsymbol{E}})+ W_{\mathrm{curv}}(  \boldsymbol \Gamma), \vspace{6pt}\label{f5}\\
W_{\mathrm{mp}}(\overline{\boldsymbol{E}} ) =\; & 
\mu\,\|\,\mathrm{dev_3\,sym}\, \overline{\boldsymbol{E}}\,\|^2\, +  \,\mu_c \, \|\,\mathrm{skew}\, \overline{\boldsymbol{E}}\,\|^2\,   + \, \dfrac{\kappa}{2}\,\big(\mathrm{tr}\,\overline{\boldsymbol{E}}\,\big)^2
\vspace{6pt}\label{f6}\\
     =\; & \mu\,\|\,\mathrm{sym}\, \overline{\boldsymbol{E}}\,\|^2\, +  \,\mu_c \, \|\,\mathrm{skew}\, \overline{\boldsymbol{E}}\,\|^2\,   + \, \dfrac{\lambda}{2}\,\big(\mathrm{tr}\,\overline{\boldsymbol{E}}\,\big)^2
     ,\vspace{6pt} \notag\\
W_{\mathrm{curv}}( \boldsymbol{\Gamma}) =\; & \mu\,L_c^2\,\Big(\,b_1\,\|\,\mathrm{dev_3\,sym}\, \boldsymbol{\Gamma}\|^2\, +  \,b_2\, \|\,\mathrm{skew}\, \boldsymbol{\Gamma}\|^2\,   + \, b_3\big(\mathrm{tr}\,\boldsymbol{\Gamma}\big)^2\,\Big) \vspace{6pt}\label{f7}\\
 =\; & \mu\,L_c^2\,\Big(\,b_1\,\|\,\mathrm{sym}\, \boldsymbol{\Gamma}\|^2\, +  \,b_2\, \|\,\mathrm{skew}\, \boldsymbol{\Gamma}\|^2\,   + \, \big(b_3- \dfrac{b_1}{3}\big)\big(\mathrm{tr}\,\boldsymbol{\Gamma}\big)^2\,\Big)\,, \notag
\end{align}
where $\mu>0$ is the shear modulus, $ \lambda $ the Lam\'e constant,  $\kappa=\frac13(3\lambda+2\mu)$ is the bulk modulus of classical isotropic elasticity, and  $\,\mu_c\ge 0$ is the so-called  \emph{Cosserat couple modulus}, $b_1\,,\, b_2 \,,\, b_3>0$ are dimensionless constitutive coefficients and the parameter $\,L_c>0\,$ introduces an internal length which is characteristic for the material.

We remark that the model is geometrically nonlinear (since the strain measures $ \overline{\boldsymbol{E}}\,,\, \boldsymbol \Gamma $ are nonlinear functions of $ \boldsymbol\varphi, \boldsymbol Q_e $), but it is physically linear in view of \eqref{f4}-\eqref{f7}. Thus, let us denote by 
$$ \underline{\boldsymbol C} = C^{ijkl}\boldsymbol g_i\otimes
\boldsymbol g_j\otimes \boldsymbol g_k\otimes \boldsymbol g_l 
 \qquad  \mathrm{and} \qquad \underline{\boldsymbol G} =  G^{ijkl}\boldsymbol g_i\otimes
 \boldsymbol g_j\otimes \boldsymbol g_k\otimes \boldsymbol g_l $$ 
the fourth-order tensors of the elastic moduli such that 
\begin{equation}\label{f8}
\begin{array}{l}
\boldsymbol Q_e^T\boldsymbol T \,=\, \underline{\boldsymbol C} : \overline{\boldsymbol{E}} \,=\, 2\mu\, \mathrm{dev_3\,sym}\, \overline{\boldsymbol{E}}  +  2\mu_c \, \mathrm{skew}\, \overline{\boldsymbol{E}}    + \, \kappa (\mathrm{tr}\,\overline{\boldsymbol{E}})\id_3 \,=\,  2\mu\, \mathrm{sym}\, \overline{\boldsymbol{E}}  +  2\mu_c \, \mathrm{skew}\, \overline{\boldsymbol{E}}    + \, \lambda (\mathrm{tr}\,\overline{\boldsymbol{E}})\id_3\,,
\vspace{6pt}\\
\boldsymbol Q_e^T\overline{\boldsymbol M} \,=\, \underline{\boldsymbol G} : \boldsymbol \Gamma \,=\, 2\mu\,L_c^2\,\Big(\,b_1\, \mathrm{dev_3\,sym}\, \boldsymbol{\Gamma}  +  \,b_2\, \mathrm{skew}\, \boldsymbol{\Gamma}    + \, b_3\big(\mathrm{tr}\,\boldsymbol{\Gamma}\big)\id_3\Big).
\end{array}
\end{equation}
By virtue of \eqref{f8}, we see that the tensor components are
\begin{equation}\label{f9}
\begin{array}{l}
C^{ijkl} =   \mu\, \big(g^{ik} g^{jl} + g^{il} g^{jk}\big) + \mu_c \big(g^{ik} g^{jl} - g^{il} g^{jk}\big) + \lambda\, g^{ij} g^{kl}\,,
\vspace{6pt}\\
G^{ijkl} = \mu\,L_c^2\,\Big(   b_1 \big(g^{ik} g^{jl} + g^{il} g^{jk}\big) + b_2 \big(g^{ik} g^{jl} - g^{il} g^{jk}\big) + 2 \big(b_3- \dfrac{b_1}{3}\big)\, g^{ij} g^{kl} \Big),
\end{array}
\end{equation}
which satisfy the major symmetries $\; C^{ijkl} = C^{klij} $ , $\; G^{ijkl} = G^{klij} \;$.
Hence, we have
\begin{equation}\label{f10}
\begin{array}{l}
W_{\mathrm{mp}}(\overline{\boldsymbol{E}} ) = \dfrac12 \big(\boldsymbol Q_e^T\boldsymbol T\big) : \overline{\boldsymbol{E}} = 
\dfrac12\; \overline{\boldsymbol{E}} : \underline{\boldsymbol C} : \overline{\boldsymbol{E}} \,,\qquad 
W_{\mathrm{curv}}(\boldsymbol{\Gamma} ) = \dfrac12 \big(\boldsymbol Q_e^T\overline{\boldsymbol M}\big) : \boldsymbol{\Gamma} = 
\dfrac12\; \boldsymbol{\Gamma} : \underline{\boldsymbol G} : \boldsymbol{\Gamma}\,.
\end{array}
\end{equation}
Under these assumptions, the deformation function $ \boldsymbol\varphi $ and microrotation tensor $ \boldsymbol Q_e  $ are the solution of the following minimization problem
\begin{equation}\label{f10,5}
I  =\dd\int_{\Omega_\xi}  W \big(\overline{\boldsymbol{E}} ,\boldsymbol{\Gamma}\big)\, \mathrm dV\quad\to \quad \textrm{\ \ min\ \   w.r.t.\ \ } ( \boldsymbol\varphi, \boldsymbol Q_e\,)\, .
\end{equation}
For the sake of simplicity, we assume here that no external body and surface loads are present. The existence of minimizers to this energy functional has been proved by the direct methods of the calculus of variations (see, e.g., \cite{Neff_Edinb06,Birsan-Neff-Ost_L56-2015}).

\section{Geometry and kinematics of three-dimensional Cosserat shells} \label{Sect3}

For a shell-like three-dimensional Cosserat body, the parametric representation $ \boldsymbol\Theta $ has the special form (see, e.g., \cite{Ciarlet00,Libai98,Pietraszkiewicz-book04})
\begin{equation}\label{f11}
\boldsymbol\Theta(\boldsymbol x )=\boldsymbol y_0(x_1,x_2)+x_3\, \boldsymbol n_0(x_1,x_2),
\end{equation}
where $\boldsymbol n_0=\dfrac{\boldsymbol y_{0,1}\times \boldsymbol y_{0,2}}{\|\boldsymbol y_{0,1}\times \boldsymbol y_{0,2}\|}\, $ is the unit normal vector  to the surface $\omega_\xi\,$, defined by the position vector $ \boldsymbol y_0(x_1,x_2) $. The parameter domain $ \Omega_h $ has the special form
$$\Omega_h=\left\{ (x_1,x_2,x_3) \,\Big|\,\, (x_1,x_2)\in\omega\subset \mathbb{R}^2 , \;\;  x_3 \in \Big( -\frac{h}{2} ,  \frac{h}{2}\, \Big)\, \right\} ,$$
where $ h $ is the thickness. Thus, $ (x_1,x_2) $ are curvilinear coordinates on the midsurface $ \omega_\xi= \boldsymbol y_0(\omega) $ and $ x_3 $ is the coordinate through the thickness of the shell-like body $ \Omega_\xi\, $.

We denote the covariant and contravariant base vectors in the tangent plane of $ \omega_\xi $ as usual by
\[ 
 \boldsymbol a_\alpha=\,\dfrac{\partial \boldsymbol y_0}{\partial x_\alpha}\,=\boldsymbol y_{0,\alpha}\,,\qquad \boldsymbol a^\beta\cdot\boldsymbol a_\alpha=\delta_\alpha^\beta\,\quad (\alpha,\beta=1,2) \qquad\mbox{and set}\qquad \boldsymbol a_3=\boldsymbol a^3=\boldsymbol n_0\,.
 \]
The surface gradient and surface divergence are then defined by
\[ 
\mathrm{Grad}_s\,\boldsymbol f\,=\, \dfrac{\partial\boldsymbol f}{\partial x_\alpha}\,\otimes \boldsymbol a^\alpha= \boldsymbol f,_{\alpha}\otimes\, \boldsymbol a^\alpha\,,\qquad \mathrm{Div}_s\,\boldsymbol T\,=\, \boldsymbol T,_{\alpha} \boldsymbol a^\alpha\,.
 \]
 We introduce the first and second fundamental tensors of the surface $ \omega_\xi\, $ by
 \begin{equation}\label{f12}
\begin{array}{l}
\boldsymbol{a}:= \text{Grad}_s\,\boldsymbol{y}_0 =\boldsymbol{a}_\alpha\otimes \boldsymbol{a}^\alpha=
a_{\alpha\beta}\boldsymbol{a}^\alpha\otimes \boldsymbol{a}^\beta= a^{\alpha\beta}\boldsymbol{a}_\alpha\otimes \boldsymbol{a}_\beta ,
\vspace{4pt}\\
\boldsymbol{b}:= -\text{Grad}_s\,\boldsymbol{n}_0=-  \boldsymbol{n}_{0,\alpha}\otimes\boldsymbol{a}^\alpha= b_{\alpha\beta}\,\boldsymbol{a}^\alpha\otimes \boldsymbol{a}^\beta=b^\alpha_\beta\,\boldsymbol{a}_\alpha\otimes \boldsymbol{a}^\beta,
\end{array}
 \end{equation}
which are symmetric. We shall also need the skew-symmetric tensor $ \boldsymbol c \,$, called the alternator tensor in the tangent plane, defined by
\begin{equation}\label{f13}
\boldsymbol c:=  \dfrac{1}{a}\,\,\epsilon_{\alpha\beta}\,\boldsymbol{a}_\alpha\otimes \boldsymbol{a}_\beta = a\,\,\epsilon_{\alpha\beta}\,\boldsymbol{a}^\alpha\otimes \boldsymbol{a}^\beta,\qquad\mbox{with}\qquad 
a:=\sqrt{\mathrm{det} \big( a_{\alpha\beta} \big)}\,>0,
\end{equation}
where $\epsilon_{\alpha\beta}\,$ is the two-dimensional alternator ($\epsilon_{12}=-\epsilon_{21}=1\,,\,\epsilon_{11}=\epsilon_{22}=0$) and $ a(x_1,x_2) $ determines the elemental area of the surface $ \omega_\xi \,$. In view of \eqref{f0,5} and \eqref{f11}, we can show that (see \cite[f. (46)]{Birsan-Neff-MMS-2019})
\begin{equation}\label{f13,5}
\boldsymbol n_0 = \boldsymbol d_3^0 = \boldsymbol Q_0 \boldsymbol e_3\,.
\end{equation}

The fundamental tensors satisfy the relation of Cayley-Hamilton type 
\begin{equation}\label{f14}
\boldsymbol b^2-2H \boldsymbol b+K\boldsymbol a=\boldsymbol 0,\qquad 2H:=\mathrm{tr}\,\boldsymbol b=b^{\alpha}_{\alpha}\,,\qquad K:=\mathrm{det}\,\boldsymbol b=\mathrm{det} \big( b^{\alpha}_{\beta} \big),
\end{equation}
where $ H $ and $ K $ are the mean curvature and the Gau\ss{} curvature of the surface $ \omega_\xi\, $, respectively. We note that $ \boldsymbol a $ plays the role of the identity tensor in the tangent plane and designate by
\begin{equation}\label{f15}
\boldsymbol b^* := -\boldsymbol b + 2H  \boldsymbol a 
\end{equation}
the cofactor of $ \boldsymbol b $ in the tangent plane, since $ \boldsymbol b\,\boldsymbol b^* = K\boldsymbol a $ in view of \eqref{f14}$ _1\, $. Let us introduce the tensors 
\begin{equation}\label{f16}
\boldsymbol \mu :=   \boldsymbol a - x_3\,\boldsymbol b,\qquad 
\boldsymbol \mu^{-1} :=  \dfrac{1}{b} (\boldsymbol a - x_3\,\boldsymbol b^*),\qquad\mbox{with}\qquad 
\boldsymbol \mu \,\boldsymbol \mu^{-1} = \boldsymbol \mu^{-1} \boldsymbol \mu = \boldsymbol a,
\end{equation}
where $ b $ is the determinant 
\begin{equation}\label{f17}
b := \mathrm{det}\,\boldsymbol \mu =  1-2H\,x_3+K\,x_3^2\, .
\end{equation}
By virtue of $ \boldsymbol g_i = \boldsymbol\Theta,_i\, $ and \eqref{f11}, \eqref{f16}, we find the relations
\begin{equation}\label{f18}
\boldsymbol g_\alpha= \boldsymbol \mu\, \boldsymbol a_\alpha\,,\qquad 
\boldsymbol g^\alpha= \boldsymbol \mu^{-1} \boldsymbol a^\alpha\,,\qquad
\boldsymbol g_3 = \boldsymbol g^3 = \boldsymbol n_0\,,
\end{equation}
which are well-known in the literature on shells. Hence, we have
\begin{equation}\label{f19}
\boldsymbol \mu = \boldsymbol g_\alpha \otimes \boldsymbol a^\alpha = \boldsymbol a^\alpha \otimes \boldsymbol g_\alpha \,,\qquad 
\boldsymbol \mu^{-1} = \boldsymbol g^\alpha \otimes \boldsymbol a_\alpha = \boldsymbol a_\alpha \otimes \boldsymbol g^\alpha \,.
\end{equation}

In the derivation of the shell model we shall employ the expansion of various functions with respect to $ x_3 \,$ about 0. Therefore, we denote the derivative of functions  with respect to $ x_3 $ with a prime, i.e. $ \,f' := \dfrac{\partial f}{\partial x_3}\, $ .\\ We can decompose the deformation gradient as follows
\begin{equation}\label{f20}
\boldsymbol F_\xi = \boldsymbol F_\xi\,\id_3  = \boldsymbol F_\xi(  \boldsymbol a + \boldsymbol n_0\otimes \boldsymbol n_0) = \boldsymbol F_\xi\, \boldsymbol a + (\boldsymbol F_\xi \boldsymbol n_0)\otimes \boldsymbol n_0\,,
\end{equation}
where
\begin{align}
& \boldsymbol F_\xi \boldsymbol n_0 = ( \boldsymbol \varphi,_i\otimes\, \boldsymbol g^i )\boldsymbol n_0 = \boldsymbol \varphi,_3 = \boldsymbol \varphi' \qquad\mbox{and}
\label{f21}\vspace{6pt}\\
& \boldsymbol F_\xi \,\boldsymbol a =  (\mathrm{Grad}_s\,\boldsymbol \varphi) \boldsymbol \mu^{-1}\,.
\label{f22}
\end{align}
To prove \eqref{f22}, we use \eqref{f18}, \eqref{f19} and write
\[  
\boldsymbol F_\xi\, \boldsymbol a = ( \boldsymbol \varphi,_i\otimes\, \boldsymbol g^i )\boldsymbol a = \boldsymbol \varphi,_\alpha\otimes\, \boldsymbol g^\alpha = (\boldsymbol\varphi,_\alpha\otimes\, \boldsymbol a^\alpha) (\boldsymbol a_\beta\otimes\, \boldsymbol g^\beta) = (\mathrm{Grad}_s\,\boldsymbol \varphi) \boldsymbol \mu^{-1}\,.
\]
Substituting \eqref{f21} and \eqref{f22} into \eqref{f20}, we get
\begin{equation}\label{f23}
\boldsymbol F_\xi = (\mathrm{Grad}_s\,\boldsymbol \varphi)\, \boldsymbol \mu^{-1} +  \boldsymbol \varphi'\otimes \boldsymbol n_0\,.
\end{equation}
We shall also need the derivatives of $ \boldsymbol F_\xi  $ with respect to $ x_3\, $. These are
\begin{equation}\label{f24}
\begin{array}{l}
\boldsymbol F_\xi' =  (\mathrm{Grad}_s\,\boldsymbol \varphi')\,\boldsymbol \mu^{-1} +
(\mathrm{Grad}_s\,\boldsymbol \varphi) \big(\boldsymbol \mu^{-1}\big)'
  +  \boldsymbol \varphi''\otimes \boldsymbol n_0\,,
\vspace{6pt}\\
\boldsymbol F_\xi'' =  (\mathrm{Grad}_s\,\boldsymbol \varphi'')\,\boldsymbol \mu^{-1} +
2(\mathrm{Grad}_s\,\boldsymbol \varphi') \big(\boldsymbol \mu^{-1}\big)'
+ (\mathrm{Grad}_s\,\boldsymbol \varphi) \big(\boldsymbol \mu^{-1}\big)''
+  \boldsymbol \varphi'''\otimes \boldsymbol n_0\,.
\end{array}
\end{equation}
Differentiating \eqref{f16} with respect to $ x_3\, $, we deduce
\begin{equation}\label{f25}
\boldsymbol \mu'=-\boldsymbol b, \qquad \boldsymbol \mu''=\boldsymbol 0,
\qquad \big(\boldsymbol \mu^{-1}\big)'=  \boldsymbol \mu^{-1} \boldsymbol b \,\boldsymbol \mu^{-1}, \qquad 
\big(\boldsymbol \mu^{-1}\big)''=  2 \boldsymbol \mu^{-1} \boldsymbol b\, \boldsymbol \mu^{-1}  \boldsymbol b\, \boldsymbol \mu^{-1}.
\end{equation}
Let us take $ x_3=0 $ in relations \eqref{f23}-\eqref{f25}. In what follows, we employ the notation $ \; \boldsymbol f_0 := \boldsymbol f_{\big| x_3=0} $ for any function $ \boldsymbol f $. Thus, we have
\begin{equation}\label{f25,5}
\boldsymbol \mu_0=\boldsymbol a, \qquad  \big(\boldsymbol \mu^{-1}\big)_0= \boldsymbol a, \qquad  \big(\boldsymbol \mu^{-1}\big)'_0= \boldsymbol b, \qquad  \big(\boldsymbol \mu^{-1}\big)''_0= 2\boldsymbol b^2
\end{equation}
and 
\begin{equation}\label{f26}
\begin{array}{l}
(\boldsymbol F_\xi)_0 =  (\mathrm{Grad}_s\,\boldsymbol \varphi)_0   +  \boldsymbol \varphi'_0\otimes \boldsymbol n_0\,,
\vspace{6pt}\\
(\boldsymbol F_\xi)_0' = (\mathrm{Grad}_s\,\boldsymbol \varphi')_0   +  (\mathrm{Grad}_s\,\boldsymbol \varphi)_0\, \boldsymbol b + \boldsymbol \varphi''_0\otimes \boldsymbol n_0\,,
\vspace{6pt}\\
(\boldsymbol F_\xi)_0'' =  (\mathrm{Grad}_s\,\boldsymbol \varphi'')_0   +  2(\mathrm{Grad}_s\,\boldsymbol \varphi')_0\, \boldsymbol b +
2(\mathrm{Grad}_s\,\boldsymbol \varphi)_0 \,\boldsymbol b^2 + \boldsymbol \varphi'''_0\otimes \boldsymbol n_0\,.
\end{array}
\end{equation}
Let us write the Taylor expansion of the deformation function $ \boldsymbol \varphi(x_1,x_2,x_3) $ with respect to $ x_3 $ in the form
\begin{equation}\label{f27}
 \boldsymbol \varphi(x_1,x_2,x_3) = \boldsymbol m(x_1, x_2) + x_3\, \boldsymbol \alpha(x_1, x_2) + \dfrac{x_3^2}{2}\,\boldsymbol \beta(x_1, x_2) +  \dfrac{x_3^3}{6}\,\boldsymbol \gamma(x_1, x_2) + \cdots \;,
\end{equation}
where 
\begin{equation}\label{f28}
\boldsymbol m = \boldsymbol \varphi_{\big| x_3=0} = \boldsymbol \varphi_0\,,  \qquad 
\boldsymbol \alpha = \boldsymbol \varphi'{}_{\big| x_3=0} = \boldsymbol \varphi'_0\,, 
\qquad 
\boldsymbol \beta = \boldsymbol \varphi''{}_{\big| x_3=0} = \boldsymbol \varphi''_0\qquad\mathrm{etc.}
\end{equation}
On the other hand, we assume that the microrotation tensor $ \boldsymbol Q_e $ does not depend on $ x_3\, $, i.e.
\begin{equation}\label{f29}
\boldsymbol Q_e(x_i) = \boldsymbol Q_e(x_1, x_2).
\end{equation}
By virtue of \eqref{f26}-\eqref{f29}, we can write the strain tensor $ \overline{\boldsymbol{E}}=\boldsymbol Q_e^T \boldsymbol F_\xi -\id_3 $ and its derivatives on the midsurface $ x_3=0\, $:
\begin{equation}\label{f30}
\begin{array}{l}
\overline{\boldsymbol{E}}_0 =\boldsymbol Q_e^T \big(\boldsymbol F_\xi\big)_0 -\id_3 = 
\boldsymbol Q_e^T \big(\mathrm{Grad}_s\,\boldsymbol m + \boldsymbol \alpha\otimes \boldsymbol n_0 \big) -\id_3 \,,
\vspace{6pt}\\
\overline{\boldsymbol{E}}_0^{\,\prime} = 
\boldsymbol Q_e^T \big(\boldsymbol F_\xi\big)'_0 =
\boldsymbol Q_e^T \big[\mathrm{Grad}_s\,\boldsymbol \alpha +
\big(\mathrm{Grad}_s\,\boldsymbol m \big)\boldsymbol b + 
 \boldsymbol \beta\otimes \boldsymbol n_0 \big]\,,
\vspace{6pt}\\
\overline{\boldsymbol{E}}_0^{\,\prime\prime} = 
\boldsymbol Q_e^T \big(\boldsymbol F_\xi\big)''_0 =
\boldsymbol Q_e^T \big[\mathrm{Grad}_s\,\boldsymbol \beta +
2\big(\mathrm{Grad}_s\,\boldsymbol\alpha \big)\boldsymbol b + 
2\big(\mathrm{Grad}_s\,\boldsymbol m \big)\boldsymbol b^2 +
\boldsymbol \gamma\otimes \boldsymbol n_0 \big]\,.
\end{array}
\end{equation}

We note that the surface $ \omega_\xi $ (characterized by $ x_3=0 $) is the midsurface of the reference shell $ \Omega_\xi\, $, while $ \boldsymbol m (x_1,x_2) $ and $ \boldsymbol Q_e(x_1, x_2) $ represent the deformation vector and microrotation tensor, respectively, for this reference midsurface $ \omega_\xi \,$. Corresponding to 
$ \boldsymbol m   $ and $ \boldsymbol Q_e  $ we introduce now the \emph{elastic shell strain tensor} $\boldsymbol E^e$ and the \emph{elastic shell bending-curvature tensor} $\boldsymbol K^e$, which are usually employed in the 6-parameter shell theory \cite{Libai98,Pietraszkiewicz-book04,Eremeyev06,Birsan-Neff-MMS-2014,Birsan-Neff-L54-2014}
\begin{equation}\label{f31}
\boldsymbol E^e  :=   \boldsymbol Q_e^T\mathrm{Grad}_s\boldsymbol m - \boldsymbol a, \qquad 
\boldsymbol K^e  :=   
\,\mathrm{axl}\big(\boldsymbol Q_e^T\boldsymbol Q_{e,\alpha}\big)\otimes \boldsymbol a^\alpha.
\end{equation}
These strain measures describe the deformation of the midsurface $ \omega_\xi\, $, see e.g. \cite{Birsan-Neff-L57-2016,Birsan-Neff-L58-2017}. With the help of \eqref{f31}$_1 $ and the decomposition $ \id_3= \boldsymbol a + \boldsymbol n_0\otimes \boldsymbol n_0 $ we can write the relation \eqref{f30}$_1 $ in the form
\begin{equation}\label{f32}
\overline{\boldsymbol{E}}_0 \;=\; \boldsymbol E^e + \big( \boldsymbol Q_e^T  \boldsymbol \alpha - \boldsymbol n_0\big)\otimes \boldsymbol n_0 \;= \;
 \boldsymbol E^e + \boldsymbol Q_e^T \big(  \boldsymbol \alpha - \boldsymbol d_3 \big)\otimes \boldsymbol n_0\,.
\end{equation}

In the same way, we can compute the wryness tensor $ \boldsymbol \Gamma $ and its derivatives on the midsurface $ x_3=0 $ in terms of the bending-curvature tensor $\boldsymbol K^e$. In view of \eqref{f18}, \eqref{f25,5} and \eqref{f29} we have
\begin{equation}\label{f33}
\begin{array}{l}
\boldsymbol \Gamma_0 = \Big(\mathrm{axl}\big(\boldsymbol Q_e^T\boldsymbol Q_{e,i}\big)\otimes \boldsymbol g^i\Big)_{x_3=0} 
= \,\mathrm{axl}\big(\boldsymbol Q_e^T\boldsymbol Q_{e,\alpha}\big)\otimes \boldsymbol a^\alpha = \boldsymbol K^e,
\vspace{6pt}\\
\boldsymbol \Gamma_0^{\,\prime} = \Big(\mathrm{axl}\big(\boldsymbol Q_e^T\boldsymbol Q_{e,i}\big)\otimes \boldsymbol g^i\Big)'_{x_3=0} 
= \,\mathrm{axl}\big(\boldsymbol Q_e^T\boldsymbol Q_{e,\alpha}\big)\otimes \big[ \big(\boldsymbol \mu^{-1}\big)'_0 \;\boldsymbol a^\alpha\big] = 
\big[\mathrm{axl}\big(\boldsymbol Q_e^T\boldsymbol Q_{e,\alpha}\big)\otimes \boldsymbol a^\alpha \big]\boldsymbol b\,
= \boldsymbol K^e\boldsymbol b,
\vspace{6pt}\\
\boldsymbol \Gamma_0^{\,\prime\prime} =  \,\mathrm{axl}\big(\boldsymbol Q_e^T\boldsymbol Q_{e,\alpha}\big)\otimes \big[ \big(\boldsymbol \mu^{-1}\big)''_0 \;\boldsymbol a^\alpha\big] = 2
\big[\mathrm{axl}\big(\boldsymbol Q_e^T\boldsymbol Q_{e,\alpha}\big)\otimes \boldsymbol a^\alpha \big]\boldsymbol b^2\,
= 2\boldsymbol K^e\boldsymbol b^2.
\end{array}
\end{equation}
These expressions will be useful in the sequel.

\section{Derivation of the two-dimensional shell model} \label{Sect4}

In order to obtain the expression of the elastically stored energy density for the two-dimensional shell model, we shall integrate the strain energy density $ W $ over the thickness and then perform some simplifications, suggested by the classical shell theory. Thus, in view of \eqref{f10,5} the total elastically stored strain-energy is
\begin{equation}\label{f34}
I  =\dd\int_{\Omega_\xi}  W \big(\overline{\boldsymbol{E}} ,\boldsymbol{\Gamma}\big)\, \mathrm dV = \int_{\omega_\xi}\Big( \int_{-h/2}^{h/2} W \big(\overline{\boldsymbol{E}} ,\boldsymbol{\Gamma}\big)\,b(x_1,x_2,x_3)\, \mathrm d x_3\Big)\mathrm d a,
\end{equation}
where $ b(x_i) $ is given by \eqref{f17} and $ \mathrm{d}a = a(x_1,x_2)\, \mathrm{d}x_1\mathrm{d}x_2 = \sqrt{\mathrm{det}(a_{\alpha\beta}) }\,\mathrm{d}x_1\mathrm{d}x_2 $ is the elemental area of the midsurface $ \omega_\xi\, $.

\subsection{Integration over the thickness} \label{Sect4.1}

With a view toward integrating with respect to $ x_3 $\,, we expand the integrand from \eqref{f34} in the form
\[  
Wb = \big(Wb\big)_0 + x_3\,\big(Wb\big)'_0 + \dfrac12\, x_3^2\,\big(Wb\big)''_0 + O(x_3^3)
\]
and find
\begin{equation}\label{f35}
\int_{-h/2}^{h/2} Wb\,\mathrm{d}x_3\, = \,h\,\big(Wb\big)_0 \,+\,  \dfrac{h^3}{24}\, \big(Wb\big)''_0 + o(h^3).
\end{equation}
By differentiating \eqref{f17} we get $ b_0=1 $ , $ b_0'=-2H $ , $ b_0''= 2K \,$. Hence, we have
\begin{equation}\label{f36}
\begin{array}{l}
\big(Wb\big)_0 = W_0\, b_0 = W_0\,,
\vspace{6pt}\\
\big(Wb\big)_0' = \big(W'b+Wb'\big)_0 = W_0'-2H\,W_0\,,
\vspace{6pt}\\
(Wb)_0'' = W_0'' -4H\,W_0'+2K\,W_0\,.
\end{array}
\end{equation}
Inserting \eqref{f36} into \eqref{f35} we obtain the expression
\begin{equation}\label{f37}
\int_{-h/2}^{h/2} Wb\,\mathrm{d}x_3\, = \Big(h+ \dfrac{h^3}{12}\,K\Big) W_0 + \dfrac{h^3}{24}\, \big(W_0'' -4H\,W_0'\big) + o(h^3).
\end{equation}
According to our constitutive assumptions \eqref{f5}-\eqref{f10}, we can write
\begin{equation}\label{f38}
\begin{array}{l}
W_0 \;= W_{\mathrm{mp}}(\overline{\boldsymbol{E}}_0)+ W_{\mathrm{curv}}(  \boldsymbol \Gamma_0) =
\dfrac12\, \overline{\boldsymbol{E}}_0 : \underline{\boldsymbol C} : \overline{\boldsymbol{E}}_0 +
\dfrac12\, \boldsymbol{\Gamma}_0 : \underline{\boldsymbol G} : \boldsymbol{\Gamma}_0\,
= \dfrac12 \big(\boldsymbol Q_e^T\boldsymbol T_0\big) : \overline{\boldsymbol{E}}_0 + \dfrac12 \big(\boldsymbol Q_e^T\overline{\boldsymbol M}_0\big) : \boldsymbol{\Gamma}_0\,,
\vspace{6pt}\\
W_0' \;= \, \overline{\boldsymbol{E}}'_0 : \underline{\boldsymbol C} : \overline{\boldsymbol{E}}_0 +
 \boldsymbol{\Gamma}'_0 : \underline{\boldsymbol G} : \boldsymbol{\Gamma}_0\,
= \big(\boldsymbol Q_e^T\boldsymbol T_0\big) : \overline{\boldsymbol{E}}'_0 + \big(\boldsymbol Q_e^T\overline{\boldsymbol M}_0\big) : \boldsymbol{\Gamma}'_0\,,
\vspace{6pt}\\
W_0''\, =  \, \overline{\boldsymbol{E}}''_0 : \underline{\boldsymbol C} : \overline{\boldsymbol{E}}_0 + \overline{\boldsymbol{E}}'_0 : \underline{\boldsymbol C} : \overline{\boldsymbol{E}}'_0 +
\boldsymbol{\Gamma}''_0 : \underline{\boldsymbol G} : \boldsymbol{\Gamma}_0 + \boldsymbol{\Gamma}'_0 : \underline{\boldsymbol G} : \boldsymbol{\Gamma}'_0
\vspace{6pt}\\
\qquad =  \, \big(\boldsymbol Q_e^T\boldsymbol T_0\big) : \overline{\boldsymbol{E}}''_0 
+ \big(\boldsymbol Q_e^T\boldsymbol T'_0\big) : \overline{\boldsymbol{E}}'_0
 + \big(\boldsymbol Q_e^T\overline{\boldsymbol M}_0\big) : \boldsymbol{\Gamma}''_0
 +  \big(\boldsymbol Q_e^T\overline{\boldsymbol M}'_0\big) : \boldsymbol{\Gamma}'_0
 \,.
\end{array}
\end{equation}
If we use the relations \eqref{f30}-\eqref{f33} in \eqref{f38} and substitute this in \eqref{f37}, we deduce the following successive expressions
\[  
\begin{array}{rl}
\displaystyle\int_{-h/2}^{h/2} Wb\,\mathrm{d}x_3\, =& \, 
\dfrac12 \Big(h+ \dfrac{h^3}{12}\,K\Big) \big[ 
\boldsymbol Q_e^T\boldsymbol T_0: \big(\boldsymbol E^e + \big( \boldsymbol Q_e^T  \boldsymbol \alpha - \boldsymbol n_0\big)\otimes \boldsymbol n_0\big) + 
\big(\boldsymbol Q_e^T\overline{\boldsymbol M}_0\big) : \boldsymbol K^e
\big]
\vspace{6pt}\\
&+ \dfrac{h^3}{24} \Big\{
\boldsymbol T_0: \big[\mathrm{Grad}_s\,\boldsymbol \beta +
2\big(\mathrm{Grad}_s\,\boldsymbol\alpha \big)\boldsymbol b + 
2\big(\mathrm{Grad}_s\,\boldsymbol m \big)\boldsymbol b^2 +
\boldsymbol \gamma\otimes \boldsymbol n_0 \big]
\vspace{6pt}\\
&+ \boldsymbol T'_0: \big[\mathrm{Grad}_s\,\boldsymbol \alpha +
\big(\mathrm{Grad}_s\,\boldsymbol m \big)\boldsymbol b + 
\boldsymbol \beta\otimes \boldsymbol n_0 \big] 
+ 2 \big(\boldsymbol Q_e^T\overline{\boldsymbol M}_0\big) : \big(\boldsymbol K^e\boldsymbol b^2\big)
+ \big(\boldsymbol Q_e^T\overline{\boldsymbol M}'_0\big) : \big(\boldsymbol K^e\boldsymbol b\big)
\vspace{6pt}\\
&-4H \,\boldsymbol T_0: \big[\mathrm{Grad}_s\,\boldsymbol \alpha +
\big(\mathrm{Grad}_s\,\boldsymbol m \big)\boldsymbol b + 
\boldsymbol \beta\otimes \boldsymbol n_0 \big]
-4H  \big(\boldsymbol Q_e^T\overline{\boldsymbol M}_0\big) : \big(\boldsymbol K^e\boldsymbol b\big)
\Big\} \; + \; o(h^3)
\end{array}
\]
or, using the decomposition $ \boldsymbol T_0 = \boldsymbol T_0 \boldsymbol a + \boldsymbol T_0 \boldsymbol n_0\otimes \boldsymbol n_0\, $,
\[  
\begin{array}{rl}
\displaystyle\int_{-h/2}^{h/2} Wb\,\mathrm{d}x_3 & = \, 
\dfrac12 \Big(h+ \dfrac{h^3}{12}\,K\Big) \big[ 
\big(\boldsymbol Q_e^T\boldsymbol T_0 \boldsymbol a\big) : \boldsymbol E^e + \big(\boldsymbol Q_e^T\boldsymbol T_0 \boldsymbol n_0\big) \cdot \big( \boldsymbol Q_e^T  \boldsymbol \alpha - \boldsymbol n_0\big) + 
\big(\boldsymbol Q_e^T\overline{\boldsymbol M}_0\big) : \boldsymbol K^e
\big]
\vspace{6pt}\\
&+ \dfrac{h^3}{24} \Big\{
\big(\boldsymbol T_0 \boldsymbol a\big) : \big[\mathrm{Grad}_s\,\boldsymbol \beta +
2\big(\mathrm{Grad}_s\,\boldsymbol\alpha \big)\boldsymbol b + 
2\big(\mathrm{Grad}_s\,\boldsymbol m \big)\boldsymbol b^2 \big]
 + \big(\boldsymbol T_0 \boldsymbol n_0\big) \cdot \boldsymbol \gamma
\vspace{6pt}\\
& + \big(\boldsymbol T'_0 \boldsymbol a\big) : \big[\mathrm{Grad}_s\,\boldsymbol \alpha +
\big(\mathrm{Grad}_s\,\boldsymbol m \big)\boldsymbol b \big] 
+ \big(\boldsymbol T'_0 \boldsymbol n_0\big) \cdot \boldsymbol \beta
+ 2 \big(\boldsymbol Q_e^T\overline{\boldsymbol M}_0\big) : \big(\boldsymbol K^e\boldsymbol b^2\big)
+ \big(\boldsymbol Q_e^T\overline{\boldsymbol M}'_0\big) : \big(\boldsymbol K^e\boldsymbol b\big)
\vspace{6pt}\\
&-4H \big(\boldsymbol T_0 \boldsymbol a\big) : \big[\mathrm{Grad}_s\,\boldsymbol \alpha +
\big(\mathrm{Grad}_s\,\boldsymbol m \big)\boldsymbol b \big] 
-4H \big(\boldsymbol T_0 \boldsymbol n_0\big) \cdot \boldsymbol \beta
-4H  \big(\boldsymbol Q_e^T\overline{\boldsymbol M}_0\big) : \big(\boldsymbol K^e\boldsymbol b\big)
\Big\} \; + \; o(h^3).
\end{array}
\]
Making some further calculations using \eqref{f14} and \eqref{f15}, we obtain
\begin{equation}\label{f39}
\begin{array}{rl}
\displaystyle\int_{-h/2}^{h/2} Wb\,\mathrm{d}x_3 & = \, 
\dfrac12 \Big(h- K\, \dfrac{h^3}{12}\Big) \big[ 
\big(\boldsymbol Q_e^T\boldsymbol T_0 \boldsymbol a\big) : \boldsymbol E^e + 
\big(\boldsymbol Q_e^T\overline{\boldsymbol M}_0\big) : \boldsymbol K^e
\big] 
+ \dfrac12 \Big(h+ \dfrac{h^3}{12}\,K\Big) \big(\boldsymbol T_0 \boldsymbol n_0\big) \cdot \big(  \boldsymbol \alpha - \boldsymbol d_3\big)
\vspace{6pt}\\
& + \dfrac{h^3}{24} \Big\{ \big(\boldsymbol T'_0 \boldsymbol a\big) : \big[\mathrm{Grad}_s\,\boldsymbol \alpha +
\big(\mathrm{Grad}_s\,\boldsymbol m \big)\boldsymbol b \big] 
+ \big(\boldsymbol T'_0 \boldsymbol n_0\big) \cdot \boldsymbol \beta
+ \big(\boldsymbol Q_e^T\overline{\boldsymbol M}'_0\big) : \big(\boldsymbol K^e\boldsymbol b\big)
\vspace{6pt}\\
&+ 
\big(\boldsymbol T_0 \boldsymbol a\big) : \big[\mathrm{Grad}_s\,\boldsymbol \beta -
2\big(\mathrm{Grad}_s\,\boldsymbol\alpha \big)\boldsymbol b^*  -
2K\big(\boldsymbol Q_e\boldsymbol a \big) \big]
+ \big(\boldsymbol T_0 \boldsymbol n_0\big) \cdot \big(\boldsymbol \gamma- 4H\boldsymbol\beta \big)
\Big\} \; + \; o(h^3).
\end{array}
\end{equation}

\subsection{Reduced form of the strain energy density} \label{Sect4.2}

The expression \eqref{f39} of the strain energy density per unit area of $ \omega_\xi $ can be further reduced, provided we make some assumptions and simplifications which are common in the classical shell theory. Thus, let us denote by $ \boldsymbol t^{\pm} $ the stress vectors on the major faces (upper and lower surfaces) of the shell, given by $ x_3=\pm\frac{h}{2} $\,. We notice that $ \boldsymbol n_0 $ is orthogonal to the major faces and write
\[  
\begin{array}{l}
 \boldsymbol t^{+} \,=\, \boldsymbol T\big(x_\alpha\,,\,\dfrac{h}{2}\;\big)\, \boldsymbol n_0 \,=\, \boldsymbol T_0 \boldsymbol n_0 + 
 \dfrac{h}{2}\,\boldsymbol T'_0 \boldsymbol n_0 +
  \dfrac{h^2}{8}\,\boldsymbol T''_0 \boldsymbol n_0 + O(h^3),
\vspace{6pt}\\
 \boldsymbol t^{-}   \,=\, \boldsymbol T\big(x_\alpha\,,\,\dfrac{-h}{2}\;\big)\, (-\boldsymbol n_0) \,=\, -\boldsymbol T_0 \boldsymbol n_0 + 
 \dfrac{h}{2}\,\boldsymbol T'_0 \boldsymbol n_0 -
 \dfrac{h^2}{8}\,\boldsymbol T''_0 \boldsymbol n_0 + O(h^3),
\end{array}
\]
which yields
\begin{equation}\label{f40}
\boldsymbol t^{+} + \boldsymbol t^{-} \,=\, h\,\boldsymbol T'_0 \boldsymbol n_0 + O(h^3)
\qquad\mathrm{and}\qquad 
\boldsymbol t^{+} - \boldsymbol t^{-} \,=\, 2\,\boldsymbol T_0 \boldsymbol n_0 + O(h^2).
\end{equation}
We assume as in the classical theory that $ \boldsymbol t^{\pm} $ are of order $ O(h^3) $ and from \eqref{f40} we find
\begin{equation}\label{f41}
\boldsymbol T_0 \boldsymbol n_0 = O(h^2)
\qquad\mathrm{and}\qquad
\boldsymbol T'_0 \boldsymbol n_0 = O(h^2).
\end{equation}
On the basis of \eqref{f41} and following the same rational as in the classical shell theory (see, e.g. \cite{Steigmann13}), we shall neglect these quantities and replace 
\begin{equation}\label{f42}
\boldsymbol T_0 \boldsymbol n_0 = \boldsymbol 0 \qquad\mathrm{and}\qquad 
\boldsymbol T_0^{\,\prime} \boldsymbol n_0 = \boldsymbol 0 
\end{equation}
in all terms of the energy density \eqref{f39}. 
Moreover, we regard the relations \eqref{f42} as two equations for the determination of the vectors $ \boldsymbol\alpha $ and $ \boldsymbol\beta $ in the expansion \eqref{f27}. Thus, from \eqref{f39} and \eqref{f42} we obtain
\begin{equation}\label{f43}
\begin{array}{rl}
\displaystyle\int_{-h/2}^{h/2} Wb\,\mathrm{d}x_3 \; = & 
\dfrac12 \Big(h- K\, \dfrac{h^3}{12}\Big) \big[ 
\big(\boldsymbol Q_e^T\boldsymbol T_0 \boldsymbol a\big) : \boldsymbol E^e + 
\big(\boldsymbol Q_e^T\overline{\boldsymbol M}_0\big) : \boldsymbol K^e
\big] 
\vspace{6pt}\\
& + \dfrac{h^3}{24} \Big\{ \big(\boldsymbol T'_0 \boldsymbol a\big) : \big[\mathrm{Grad}_s\,\boldsymbol \alpha +
\big(\mathrm{Grad}_s\,\boldsymbol m \big)\boldsymbol b \big] 
+ \big(\boldsymbol Q_e^T\overline{\boldsymbol M}'_0\big) : \big(\boldsymbol K^e\boldsymbol b\big)
\vspace{6pt}\\
& \qquad \;\; + 
\big(\boldsymbol T_0 \boldsymbol a\big) : \big[\mathrm{Grad}_s\,\boldsymbol \beta -
2\big(\mathrm{Grad}_s\,\boldsymbol\alpha \big)\boldsymbol b^*  -
2K\big(\boldsymbol Q_e\boldsymbol a \big) \big]
\Big\}.
\end{array}
\end{equation}
In view of \eqref{f30}-\eqref{f32}, the equations \eqref{f42} can be written in the form
\begin{equation}\label{f44}
\begin{array}{l}
 \Big[\, \underline{\boldsymbol C} :
\big(\boldsymbol E^e + \big( \boldsymbol Q_e^T  \boldsymbol \alpha - \boldsymbol n_0\big)\otimes \boldsymbol n_0\big)\Big] \boldsymbol n_0 =\boldsymbol 0,
\vspace{6pt}\\
 \Big[\, \underline{\boldsymbol C} :
 \big(\boldsymbol Q_e^T \mathrm{Grad}_s\,\boldsymbol \alpha+ (\boldsymbol E^e+ \boldsymbol a)\boldsymbol b + \boldsymbol Q_e^T  \boldsymbol \beta \otimes \boldsymbol n_0\big)\Big] \boldsymbol n_0 =\boldsymbol 0.
\end{array}
\end{equation}
The first equation \eqref{f44}$ _1 $ can be used to determine the vector $ \boldsymbol\alpha $\,: we obtain successively
\[  
\Big[ (\mu+\mu_c)\boldsymbol a + (\lambda+2\mu)\, \boldsymbol n_0\otimes \boldsymbol n_0\Big] \big( \boldsymbol Q_e^T  \boldsymbol \alpha - \boldsymbol n_0\big) = - \big( \underline{\boldsymbol C} :
\boldsymbol E^e\big)  \boldsymbol n_0\,,
\]
or equivalently,
\[  
\boldsymbol Q_e^T  \boldsymbol \alpha - \boldsymbol n_0 =
- \Big[ \dfrac{1}{\mu+\mu_c}\;\boldsymbol a + \dfrac{1}{\lambda+2\mu}\; \boldsymbol n_0\otimes \boldsymbol n_0\Big]
\big[
(\mu-\mu_c)\big( \boldsymbol n_0 \boldsymbol E^e\big) 
+\lambda \big( \mathrm{tr}\, \boldsymbol E^e\big)\boldsymbol n_0
\big],
\]
which yields (since $ \boldsymbol Q_e\boldsymbol n_0 = \boldsymbol Q_e\boldsymbol d_3^0 =\boldsymbol d_3$)
\begin{equation}\label{f45}
\boldsymbol\alpha = \Big(1- \dfrac{\lambda}{\lambda+2\mu}\,\mathrm{tr}\, \boldsymbol E^e\Big)\boldsymbol d_3 \;-\; 
\dfrac{\mu-\mu_c}{\mu+\mu_c}\; \boldsymbol Q_e\big( \boldsymbol n_0 \boldsymbol E^e\big).
\end{equation}
Further, we solve the second equation \eqref{f44}$ _2 $ to determine the vector $ \boldsymbol\beta $. To this aim, we insert $ \boldsymbol\alpha $ given by \eqref{f45}
into \eqref{f44}$ _2 $ and (in order to avoid quadratic terms and derivatives of the strain measures $ \boldsymbol E^e, \boldsymbol K^e $) we use the approximation
\[
\boldsymbol Q_e^T \mathrm{Grad}_s\boldsymbol\alpha\; \simeq\; 
\boldsymbol Q_e^T \mathrm{Grad}_s\boldsymbol d_3\,.
\]
Since $ \;\boldsymbol Q_e^T \mathrm{Grad}_s\boldsymbol d_3 = \boldsymbol c \boldsymbol K^e -\boldsymbol b \;$ (see \cite[f. (70)]{Birsan-Neff-MMS-2019}), we use
\begin{equation}\label{f46}  
\boldsymbol Q_e^T \mathrm{Grad}_s\boldsymbol\alpha\;= \;  \boldsymbol c \boldsymbol K^e -\boldsymbol b 
\end{equation}
and the equation \eqref{f44}$ _2 $ becomes
\[  
 \Big[\, \underline{\boldsymbol C} :
 \big(\boldsymbol E^e\boldsymbol b + \boldsymbol c \boldsymbol K^e+ \boldsymbol Q_e^T  \boldsymbol \beta \otimes \boldsymbol n_0\big)\Big] \boldsymbol n_0 =\boldsymbol 0,
\]
which can be solved similarly as the equation \eqref{f44}$ _1 $ and yields
\begin{equation}\label{f47}
\boldsymbol \beta = 
- \dfrac{\lambda}{\lambda+2\mu}\,\mathrm{tr}\,\big(\boldsymbol E^e\boldsymbol b + \boldsymbol c \boldsymbol K^e \big)\, \boldsymbol d_3
 \;-\; 
 \dfrac{\mu-\mu_c}{\mu+\mu_c}\; \boldsymbol Q_e\big( \boldsymbol n_0 \boldsymbol E^e\boldsymbol b\big).
\end{equation}
In view of \eqref{f45}-\eqref{f47}, we can write the tensors $ \overline{\boldsymbol{E}}_0 $ and $
\overline{\boldsymbol{E}}_0^{\,\prime}  $ in \eqref{f32} and \eqref{f30}$ _2 $ in compact form
\begin{equation}\label{f48}
\begin{array}{l}
\overline{\boldsymbol{E}}_0 = \boldsymbol E^e -\Big[\;
\dfrac{\lambda}{\lambda+2\mu}\,\big(\mathrm{tr}\,\boldsymbol E^e\big)\,\boldsymbol n_0
+ \dfrac{\mu-\mu_c}{\mu+\mu_c}\,\big( \boldsymbol n_0 \boldsymbol E^e\big) 
\Big] \otimes \boldsymbol n_0 = L_{n_0} \big( \boldsymbol E^e\big),
\vspace{6pt}\\
\overline{\boldsymbol{E}}_0^{\,\prime} = 
\big(\boldsymbol E^e\boldsymbol b + \boldsymbol c \boldsymbol K^e \big) -\Big[\;
\dfrac{\lambda}{\lambda+2\mu}\,\mathrm{tr}\big(\boldsymbol E^e\boldsymbol b + \boldsymbol c \boldsymbol K^e \big)\,\boldsymbol n_0
+ \dfrac{\mu-\mu_c}{\mu+\mu_c}\,\big( \boldsymbol n_0 \boldsymbol E^e\boldsymbol b\big)
\Big] \otimes \boldsymbol n_0 = L_{n_0} \big(\boldsymbol E^e\boldsymbol b + \boldsymbol c \boldsymbol K^e \big),
\end{array}
\end{equation}
where we have denoted for convenience with $ L_{n_0} $ the following linear operator
\begin{equation}\label{f49}
 L_{n_0}(\boldsymbol X) := \,
 \boldsymbol X \,-\,
  \,\dfrac{\lambda}{\lambda+2\mu}\,\big(\mathrm{tr}\,\boldsymbol X\big)\,\boldsymbol n_0 \otimes \boldsymbol n_0
  \,-\,
 \dfrac{\mu-\mu_c}{\mu+\mu_c}\,\big( \boldsymbol n_0 \boldsymbol X\big)
 \otimes \boldsymbol n_0
   \qquad \text{for any}\qquad\boldsymbol X = X_{i\alpha}\boldsymbol a^i\otimes \boldsymbol a^\alpha .
\end{equation}
To write the strain-energy density in a condensed form, we designate by 
\begin{equation}\label{f49,5}
\begin{array}{rl}
W_{\mathrm{mixt}}(\boldsymbol X,\boldsymbol Y) :=&   \mu\, (\mathrm{sym}\, \boldsymbol X) : (\mathrm{sym}\, \boldsymbol Y) +  \mu_c (\mathrm{skew}\, \boldsymbol X) : (\mathrm{skew}\, \boldsymbol Y) +\,\dfrac{\lambda\,\mu}{\lambda+2\mu}\,\big( \mathrm{tr} \boldsymbol X\big)\,\big(\mathrm{tr} \boldsymbol Y\big) 
\vspace{6pt}\\
= &   \mu\,(\mathrm{dev_3\,sym}\, \boldsymbol X) : (\mathrm{dev_3\,sym}\, \boldsymbol Y)  +  \mu_c  (\mathrm{skew} \boldsymbol X) : (\mathrm{skew}\, \boldsymbol Y)  +\,\dfrac{2\mu(2\lambda+\mu)}{3(\lambda+2\mu)}\,\big( \mathrm{tr} \boldsymbol X\big)\,\big(\mathrm{tr} \boldsymbol Y\big)
\end{array}
\end{equation}
the bilinear form corresponding to the quadratic form
\begin{equation}\label{f49,6}
\begin{array}{rl}
W_{\mathrm{mixt}}(\boldsymbol X) :=&   W_{\mathrm{mixt}}(\boldsymbol X,\boldsymbol X) \,\,=\,\,
 W_{\mathrm{mp}}(\boldsymbol X) - \, \dfrac{\lambda^2}{2(\lambda+2\mu)}\,\big( \mathrm{tr} \boldsymbol X\big)^2
\vspace{6pt}\\
= &  \mu\,\|\, \mathrm{sym}\, \boldsymbol X\,\|^2 +  \mu_c\| \,\mathrm{skew}\, \boldsymbol X\,\|^2 +\,\dfrac{\lambda\,\mu}{\lambda+2\mu}\,\big( \mathrm{tr} \boldsymbol X\big)^2 .
\end{array}
\end{equation}
For Cosserat shells, it is convenient to introduce the following bilinear form 
\begin{equation}\label{f50}
W_{\mathrm{Coss}}(\boldsymbol X,\boldsymbol Y) :=   W_{\mathrm{mixt}}(\boldsymbol X,\boldsymbol Y) - \, \dfrac{(\mu-\mu_c)^2}{2(\mu+\mu_c)}\,\big(\boldsymbol n_0  \boldsymbol X\big)\cdot\big(\boldsymbol n_0  \boldsymbol Y\big)
\end{equation}
for any two tensors of the form $ \boldsymbol X = X_{i\alpha}\boldsymbol a^i\otimes \boldsymbol a^\alpha $, $\; \boldsymbol Y = Y_{i\alpha}\boldsymbol a^i\otimes \boldsymbol a^\alpha $, and the corresponding quadratic form 
\begin{equation}\label{f51}
W_{\mathrm{Coss}}(\boldsymbol X) :=   W_{\mathrm{Coss}}(\boldsymbol X,\boldsymbol X) \;=\; W_{\mathrm{mixt}}(\boldsymbol X) - \, \dfrac{(\mu-\mu_c)^2}{2(\mu+\mu_c)}\,\|\boldsymbol n_0  \boldsymbol X\|^2,
\end{equation}
where $ \boldsymbol n_0  \boldsymbol X = X_{3\alpha} \boldsymbol a^\alpha$. We shall prove later that the quadratic form 
$
W_{\mathrm{Coss}}(\boldsymbol X) $ is positive definite, see \eqref{f99}. 

With these notations, we can prove by a straightforward calculation the following useful relation 
\begin{equation}\label{f52}
W_{\mathrm{Coss}}(\boldsymbol X)  
\;=\; \dfrac12\;
\boldsymbol X : \underline{\boldsymbol C} : L_{n_0}(\boldsymbol X)  
 \qquad \text{for any}\quad\boldsymbol X = X_{i\alpha}\boldsymbol a^i\otimes \boldsymbol a^\alpha .
\end{equation}
Indeed, we have from \eqref{f10}, \eqref{f49},  \eqref{f49,6}, \eqref{f51}
\[  
\begin{array}{rl}
\boldsymbol X : \underline{\boldsymbol C} : L_{n_0}(\boldsymbol X) &= \boldsymbol X : \underline{\boldsymbol C} : \boldsymbol X - 
\boldsymbol X :\underline{\boldsymbol C} : \Big[
\,\dfrac{\lambda}{\lambda+2\mu}\,\big(\mathrm{tr}\,\boldsymbol X\big)\,\boldsymbol n_0 \otimes \boldsymbol n_0
+
\dfrac{\mu-\mu_c}{\mu+\mu_c}\,\big( \boldsymbol n_0 \boldsymbol X\big)
\otimes \boldsymbol n_0\Big]
\vspace{6pt}\\
&=2W_{\mathrm{mp}}(\boldsymbol X) - \boldsymbol X : \Big[ 
\,\dfrac{\lambda^2}{\lambda+2\mu}\,\big(\mathrm{tr}\,\boldsymbol X\big)\,\id_3
+
(\mu-\mu_c)\,\big( \boldsymbol n_0 \boldsymbol X\big)
\otimes \boldsymbol n_0
+
\dfrac{(\mu-\mu_c)^2}{\mu+\mu_c}\,\boldsymbol n_0 \otimes \big( \boldsymbol n_0 \boldsymbol X\big)
 \Big]
 \vspace{6pt}\\
 &=2W_{\mathrm{mp}}(\boldsymbol X) -  
 \,\dfrac{\lambda^2}{\lambda+2\mu}\,\big(\mathrm{tr}\,\boldsymbol X\big)^2
- 
 \dfrac{(\mu-\mu_c)^2}{\mu+\mu_c}\,\| \boldsymbol n_0 \boldsymbol X\|^2
  \vspace{6pt}\\
  &=2W_{\mathrm{mixt}}(\boldsymbol X) -  
  \dfrac{(\mu-\mu_c)^2}{\mu+\mu_c}\,\| \boldsymbol n_0 \boldsymbol X\|^2 \;=\; 2W_{\mathrm{Coss}}(\boldsymbol X)
\end{array}
\]
and the relation \eqref{f52} is proved. 

Now, we can simplify the terms appearing in the strain-energy density \eqref{f43}: making use of \eqref{f42}, \eqref{f46}, \eqref{f48} and \eqref{f52} we find
\begin{equation}\label{f53}
\big(\boldsymbol Q_e^T\boldsymbol T_0 \boldsymbol a\big) : \boldsymbol E^e = 
\boldsymbol E^e : \big(\boldsymbol Q_e^T\boldsymbol T_0 \big) = \boldsymbol E^e : \big(\underline{\boldsymbol C} : 
\overline{\boldsymbol{E}}_0\big) = \boldsymbol E^e : \underline{\boldsymbol C} : 
L_{n_0} \big( \boldsymbol E^e\big) = 2W_{\mathrm{Coss}}\big( \boldsymbol E^e\big)
\end{equation}
and
\begin{equation}\label{f54}
\begin{array}{l}
\big(\boldsymbol T'_0 \boldsymbol a\big) : \big[\mathrm{Grad}_s\,\boldsymbol \alpha +
\big(\mathrm{Grad}_s\,\boldsymbol m \big)\boldsymbol b \big] = 
\big(\boldsymbol Q_e^T\boldsymbol T'_0 \boldsymbol a\big) : \big[ \big( \boldsymbol c \boldsymbol K^e -\boldsymbol b \big) + \big( \boldsymbol E^e + \boldsymbol a\big)\boldsymbol b\big] 
 = \big(\boldsymbol Q_e^T\boldsymbol T'_0 \big) : \big(\boldsymbol E^e\boldsymbol b + \boldsymbol c \boldsymbol K^e \big) 
 \vspace{6pt}\\
 \qquad\qquad
= \big(\boldsymbol E^e\boldsymbol b + \boldsymbol c \boldsymbol K^e \big) :  
 \big(\underline{\boldsymbol C} : \overline{\boldsymbol{E}}_0^{\,\prime}\big) =
 \big(\boldsymbol E^e\boldsymbol b + \boldsymbol c \boldsymbol K^e \big) :  
 \underline{\boldsymbol C} : L_{n_0}\big(\boldsymbol E^e\boldsymbol b + \boldsymbol c \boldsymbol K^e \big) =
  2 W_{\mathrm{Coss}}\big(\boldsymbol E^e\boldsymbol b + \boldsymbol c \boldsymbol K^e \big)
\end{array}
\end{equation}
and
\begin{equation}\label{f55}
\begin{array}{l}
\big(\boldsymbol T_0 \boldsymbol a\big) : \big[\big(\mathrm{Grad}_s\,\boldsymbol\alpha \big)\boldsymbol b^*  
+K\big(\boldsymbol Q_e\boldsymbol a \big) \big] = 
\big(\boldsymbol Q_e^T\boldsymbol T_0 \boldsymbol a\big) : \big[ \big( \boldsymbol c \boldsymbol K^e -\boldsymbol b \big)\boldsymbol b^* + K \boldsymbol a\big] 
= 
\big(\boldsymbol Q_e^T\boldsymbol T_0 \big) :  \big( \boldsymbol c \boldsymbol K^e \boldsymbol b^* \big) 
 \vspace{6pt}\\
 \qquad\qquad
=   
\big(\underline{\boldsymbol C} : \overline{\boldsymbol{E}}_0\big) : \big( \boldsymbol c \boldsymbol K^e \boldsymbol b^* \big)   =
\big[ 2\mu\, \mathrm{sym}\, \overline{\boldsymbol{E}}_0  +  2\mu_c \, \mathrm{skew}\, \overline{\boldsymbol{E}}_0    + \, \lambda (\mathrm{tr}\,\overline{\boldsymbol{E}}_0)\id_3\big]: \big( \boldsymbol c \boldsymbol K^e \boldsymbol b^* \big)   
 \vspace{6pt}\\
 \qquad\qquad
 =   2\mu\, \mathrm{sym}\big( \boldsymbol E^e\big) : \mathrm{sym}\, \big( \boldsymbol c \boldsymbol K^e \boldsymbol b^* \big)    +  2\mu_c\, \mathrm{skew}\big( \boldsymbol E^e ) : \mathrm{skew}\, \big( \boldsymbol c \boldsymbol K^e \boldsymbol b^* \big)    +\,\dfrac{2\lambda\,\mu}{\lambda+2\mu}\, \mathrm{tr} \big(\boldsymbol E^e \big)\,\mathrm{tr} \big( \boldsymbol c \boldsymbol K^e \boldsymbol b^* \big)    
  \vspace{6pt}\\
  \qquad\qquad
  =   2 W_{\mathrm{Coss}}\big(\boldsymbol E^e\,,\, \boldsymbol c \boldsymbol K^e \boldsymbol b^*  \big),
\end{array}
\end{equation}
since $ \,\,\mathrm{tr}\,\overline{\boldsymbol{E}}_0 = \,\dfrac{2\,\mu}{\lambda+2\mu}\,\mathrm{tr} \,\boldsymbol E^e \,\,$ and  the tensor $ \,\boldsymbol c \boldsymbol K^e \boldsymbol b^*  \,$ is a planar tensor with basis $ \{\boldsymbol a^\alpha\otimes\boldsymbol a^\beta  \} $.

Further, the two terms involving the bending-curvature tensor $ \boldsymbol K^e $ in the strain-energy density \eqref{f43} can be transformed as follows: by virtue of \eqref{f8}, \eqref{f10} and \eqref{f33} we have
\begin{equation}\label{f56}
\big(\boldsymbol Q_e^T\overline{\boldsymbol M}_0\big) : \boldsymbol K^e = \boldsymbol K^e : \big( \underline{\boldsymbol G} : \boldsymbol{\Gamma}_0 \big) 
 = \boldsymbol K^e : \underline{\boldsymbol G} : \boldsymbol K^e = 2 W_{\mathrm{curv}} \big(\boldsymbol K^e\big)
\end{equation}
and 
\begin{equation}\label{f57}
\big(\boldsymbol Q_e^T\overline{\boldsymbol M}'_0\big) : \big(\boldsymbol K^e \boldsymbol b\big) = \big(\boldsymbol K^e \boldsymbol b\big) : \big( \underline{\boldsymbol G} : \boldsymbol{\Gamma}'_0 \big) 
= \big(\boldsymbol K^e \boldsymbol b\big) : \underline{\boldsymbol G} : \big(\boldsymbol K^e \boldsymbol b\big) = 2 W_{\mathrm{curv}} \big(\boldsymbol K^e \boldsymbol b\big).
\end{equation}
Finally, the term $ (\boldsymbol T_0\boldsymbol a) : \mathrm{Grad}_s\boldsymbol\beta \; $ appearing in the strain-energy density \eqref{f43} can be discarded. To justify this, we proceed as in the classical shell theory, see e.g. \cite{Steigmann12,Steigmann13}: the three-dimensional equilibrium equation 
$\; \mathrm{Div}\,\boldsymbol T = \boldsymbol 0\; $
can be written as $\; \boldsymbol T,_i\boldsymbol g^i  = \boldsymbol 0\, $, or equivalently
\[  
\boldsymbol T,_\alpha\boldsymbol g^\alpha + \boldsymbol T^{\,\prime}\boldsymbol n_0 = \boldsymbol 0.
\]
Therefore, on the midsurface $ x_3=0 $ we have
\begin{equation}\label{f58}
\boldsymbol T_{0,\alpha}\boldsymbol a^\alpha + \boldsymbol T_0^{\,\prime}\boldsymbol n_0 = \boldsymbol 0.
\end{equation}
On the other hand,  we see that 
\[  
\boldsymbol T_{0,\alpha}\boldsymbol a^\alpha = \big( \boldsymbol T_0 \boldsymbol a + \boldsymbol T_0 \boldsymbol n_0\otimes \boldsymbol n_0\big),_\alpha \boldsymbol a^\alpha  =
 \big( \boldsymbol T_0 \boldsymbol a \big),_\alpha \boldsymbol a^\alpha
+ \boldsymbol T_0 \boldsymbol n_0\big(\boldsymbol n_{0,\alpha} \cdot \,\boldsymbol a^\alpha\big)
= \mathrm{Div}_s(\boldsymbol T_0\boldsymbol a) -2H \,\boldsymbol T_0 \boldsymbol n_0\,.
\]
Inserting the last relation into \eqref{f58} we find
\begin{equation}\label{f59}
\mathrm{Div}_s(\boldsymbol T_0\boldsymbol a) + \boldsymbol T_0^{\,\prime}\boldsymbol n_0 -2H \,\boldsymbol T_0 \boldsymbol n_0 = \boldsymbol 0.
\end{equation}
With help of \eqref{f42}, \eqref{f59} and the divergence theorem for surfaces we get
\begin{equation}\label{f60}
\begin{array}{l}
\displaystyle\int_{\omega_\xi}  \big(\boldsymbol T_0 \boldsymbol a\big) : \big(\mathrm{Grad}_s\,\boldsymbol \beta \big) \mathrm{d}a =
\int_{\omega_\xi}  \big[ \mathrm{Div}_s\big(\boldsymbol \beta(\boldsymbol T_0\boldsymbol a)\big)  - \boldsymbol \beta\cdot \mathrm{Div}_s(\boldsymbol T_0\boldsymbol a)\big] \mathrm{d}a =
 \vspace{6pt}\\
 \qquad\qquad\displaystyle
= \int_{\partial\omega_\xi} \boldsymbol \beta\big(\boldsymbol T_0\boldsymbol a\big)  \cdot \boldsymbol \nu\, \mathrm{d}\ell -
 \int_{\omega_\xi}   \boldsymbol \beta\cdot \big( 2H \,\boldsymbol T_0 \boldsymbol n_0  - \boldsymbol T_0^{\,\prime}\boldsymbol n_0 \big) \mathrm{d}a 
 = \int_{\partial\omega_\xi}  \boldsymbol \beta\cdot \big(\boldsymbol T_0\boldsymbol a\big)  \boldsymbol \nu\, \mathrm{d}\ell\,,
\end{array}
\end{equation}
where $ \boldsymbol\nu  $ is the unit normal to the boundary curve  $ \partial \omega_\xi $ lying in the tangent plane. The last integral in 
\eqref{f60} represents a prescribed constant (determined by the boundary data on $ \partial \omega_\xi $), which can be omitted, since its variation vanishes identically and thus does not influence the minimizers of the energy functional.

In conclusion, using the results \eqref{f53}-\eqref{f57} in the equation \eqref{f43} we obtain the following expression of the areal strain-energy density for Cosserat shells
\begin{equation}\label{f61}
\begin{array}{rl}
W_{\mathrm{shell}}(\boldsymbol E^e,\boldsymbol K^e) = &
\Big(h- K\, \dfrac{h^3}{12}\Big) \big[  W_{\mathrm{Coss}}\big( \boldsymbol E^e\big) + W_{\mathrm{curv}} \big(\boldsymbol K^e\big) \big]
 \vspace{6pt}\\
& + \;\dfrac{h^3}{12}\,\big[  W_{\mathrm{Coss}}\big(\boldsymbol E^e\boldsymbol b + \boldsymbol c \boldsymbol K^e \big) 
-2 W_{\mathrm{Coss}}\big(\boldsymbol E^e\,,\, \boldsymbol c \boldsymbol K^e \boldsymbol b^*  \big)
+ W_{\mathrm{curv}} \big(\boldsymbol K^e \boldsymbol b\big) \big],
\end{array}
\end{equation}
where $ W_{\mathrm{Coss}} $ is defined by \eqref{f50}, \eqref{f51} (see also equations \eqref{f88} and \eqref{f99})
and $ W_{\mathrm{curv}} $ is given in \eqref{f7}. This is the elastically stored strain-energy density for our model, which determines the constitutive equations. In Section \ref{Sect5} we shall present a useful alternative form of the 
energy $ W_{\mathrm{shell}}(\boldsymbol E^e,\boldsymbol K^e) $, together with explicit stress-strain relations (see \eqref{f100}, \eqref{f107}).

\subsection{The field equations for Cosserat shells} \label{Sect4.3}

For the sake of completeness, we record here the governing field equations of the derived shell model. 

We deduce the form of the equilibrium equations for Cosserat shells from the condition that the solution is a stationary point of the energy functional $ I $\,, i.e. we impose that the variation of the energy functional is zero:
\begin{equation}\label{f62}
\delta I =0\,, \qquad \mathrm{with}\quad I=\int_{\omega_\xi}
W_{\mathrm{shell}}(\boldsymbol E^e,\boldsymbol K^e) \,\mathrm{d}a.
\end{equation}
For simplicity we have assumed in \eqref{f62} that the external body loads are vanishing and the boundary conditions are null. To compute the variation $ \delta I  $ we write
\begin{equation}\label{f63}
\delta\, 
W_{\mathrm{shell}}(\boldsymbol E^e,\boldsymbol K^e) \, = \,\dfrac{\partial\, W_{\mathrm{shell}}}{\partial \boldsymbol E^e}\,  
: \big(\delta\boldsymbol E^e\big) +
\dfrac{\partial\, W_{\mathrm{shell}}}{\partial \boldsymbol K^e}\,: \big(\delta\boldsymbol K^e\big) 
\, = \,\big(\boldsymbol Q_e^T\boldsymbol N \big) 
: \big(\delta\boldsymbol E^e\big) +
\big(\boldsymbol Q_e^T\boldsymbol M\big) : \big(\delta\boldsymbol K^e\big) ,
\end{equation}
where we have introduced the tensors $ \boldsymbol N $ and $ \boldsymbol M $ such that
\begin{equation}\label{f64}
\boldsymbol Q_e^T\boldsymbol N = \dfrac{\partial\, W_{\mathrm{shell}}}{\partial \boldsymbol E^e}\qquad 
\textrm{and}\qquad
\boldsymbol Q_e^T\boldsymbol M = \dfrac{\partial\, W_{\mathrm{shell}}}{\partial \boldsymbol K^e}\,. 
\end{equation}
Let us denote by 
\begin{equation}\label{f65}
\boldsymbol F_s :=\, \mathrm{Grad}_s\boldsymbol m \,=\, \boldsymbol m,_\alpha \otimes\, \boldsymbol a^\alpha
\end{equation}
the \textit{shell deformation gradient} (i.e., the surface gradient of the midsurface deformation $ \boldsymbol m $). Then, in view of \eqref{f31}$ _1 $ we have $\, \boldsymbol E^e = \boldsymbol Q_e^T \boldsymbol F_s - \boldsymbol a\,$ and, hence,
\begin{equation}\label{f66}
\delta\boldsymbol E^e = \delta\big( \boldsymbol Q_e^T \boldsymbol F_s - \boldsymbol a \big) = \delta \big( \boldsymbol Q_e^T \mathrm{Grad}_s\boldsymbol m \big)
= ( \delta \boldsymbol Q_e)^T \mathrm{Grad}_s\boldsymbol m  + 
  \boldsymbol Q_e^T \mathrm{Grad}_s\big( \delta\boldsymbol m \big).
\end{equation}
To compute $ \delta\boldsymbol Q_e\, $, we notice that  the tensor $ (\delta\boldsymbol Q_e)\boldsymbol Q_e^T $ is skew-symmetric and we denote
\begin{equation}\label{f67}
\boldsymbol\Omega := (\delta\boldsymbol Q_e)\boldsymbol Q_e^T, \qquad \boldsymbol\omega:=\mathrm{axl} (\boldsymbol\Omega),\qquad
\mathrm{with}\quad \boldsymbol\Omega=\boldsymbol\omega\times \id_3\,.
\end{equation}
In the above relations, the axial vector $ \boldsymbol\omega $ is the virtual rotation vector and $ \delta\boldsymbol m$ is the virtual translation. From \eqref{f67} we get 
\begin{equation}\label{f68}
\delta\boldsymbol Q_e = \boldsymbol\Omega \,\boldsymbol Q_e
=-(\boldsymbol Q_e^T \boldsymbol\Omega)^T
\end{equation}
and substituting into \eqref{f66} we obtain
\begin{equation}\label{f69}
\delta\boldsymbol E^e = \boldsymbol Q_e^T\big( \mathrm{Grad}_s (\delta \boldsymbol m) -  \boldsymbol\Omega \,\boldsymbol F_s \big)
.
\end{equation}
Further, in order to compute $  \delta\boldsymbol K^e\, $, we recall the formula (see \cite[f. (63)]{Birsan-Neff-L57-2016})
\begin{equation}\label{f70}
\boldsymbol K^e = \dfrac12\,\big[ \boldsymbol Q_e^T \big(   \boldsymbol d_i\times \mathrm{Grad}_s\, \boldsymbol d_i\big) 
-  \boldsymbol d^0_i\times \mathrm{Grad}_s\, \boldsymbol d^0_i \,\big]
\end{equation}
and write (in view of \eqref{f68}) 
\begin{equation}\label{f71}
\delta\boldsymbol d_i = \delta \big( \boldsymbol Q_e   \boldsymbol d_i^0\big) = (\delta\boldsymbol Q_e)  \boldsymbol d_i^0 = \boldsymbol\Omega \boldsymbol Q_e  \boldsymbol d_i^0 = \boldsymbol\Omega   \boldsymbol d_i = \boldsymbol\omega\times
\boldsymbol d_i\,.
\end{equation}
Then, from \eqref{f70} it follows
\begin{equation}\label{f72}
\begin{array}{rl}
 \delta\boldsymbol K^e =&  \dfrac12 \, \delta \big[ \boldsymbol Q_e^T \big(   \boldsymbol d_i\times \mathrm{Grad}_s\, \boldsymbol d_i\big) \big]
 \vspace{6pt}\\
 = & 
 \dfrac12 \,  \big[ (\delta\boldsymbol Q_e)^T \big(   \boldsymbol d_i\times \mathrm{Grad}_s\, \boldsymbol d_i\big) 
 + \boldsymbol Q_e^T \big(   (\delta\boldsymbol d_i)\times \mathrm{Grad}_s\, \boldsymbol d_i\big)
 + \boldsymbol Q_e^T \big(   \boldsymbol d_i\times \mathrm{Grad}_s (\delta\boldsymbol d_i)\big)\big] 
 \vspace{6pt}\\
 = & 
 \dfrac12 \, \boldsymbol Q_e^T \big[ -\boldsymbol \Omega \big(   \boldsymbol d_i\times \mathrm{Grad}_s\, \boldsymbol d_i\big) 
+   ( \boldsymbol \Omega\boldsymbol d_i)\times \mathrm{Grad}_s\, \boldsymbol d_i
 +    \boldsymbol d_i\times \mathrm{Grad}_s (\boldsymbol \Omega\boldsymbol d_i)\big] 
\vspace{6pt}\\
= &  \dfrac12 \, \boldsymbol Q_e^T \big[ -\boldsymbol\omega\times \big(   \boldsymbol d_i\times \mathrm{Grad}_s\, \boldsymbol d_i\big) 
+   ( \boldsymbol\omega\times\boldsymbol d_i)\times \mathrm{Grad}_s\, \boldsymbol d_i
+    \boldsymbol d_i\times \mathrm{Grad}_s (\boldsymbol\omega\times\boldsymbol d_i)\big] .
\end{array}
\end{equation}
By virtue of the Jacobi identity for the cross product, we have
\[  
 -\boldsymbol\omega\times \big(   \boldsymbol d_i\times \mathrm{Grad}_s\, \boldsymbol d_i\big) 
 +   ( \boldsymbol\omega\times\boldsymbol d_i)\times \mathrm{Grad}_s\, \boldsymbol d_i = 
  - \boldsymbol d_i\times\big(    \boldsymbol\omega\times\mathrm{Grad}_s\, \boldsymbol d_i\big) 
\]
and inserting this in \eqref{f72} we get
\begin{equation}\label{f72,5}
\begin{array}{rl}
\delta\boldsymbol K^e =&   \dfrac12 \, \boldsymbol Q_e^T \big[ \boldsymbol d_i \times \big( \mathrm{Grad}_s (\boldsymbol\omega\times\boldsymbol d_i) -  \boldsymbol\omega\times\mathrm{Grad}_s\, \boldsymbol d_i \big)\big] .
\end{array}
\end{equation}
For the square brackets in \eqref{f72,5} we can write
\begin{equation}\label{f73}
\boldsymbol d_i \times \big( \mathrm{Grad}_s (\boldsymbol\omega\times\boldsymbol d_i) -  \boldsymbol\omega\times\mathrm{Grad}_s\, \boldsymbol d_i \big)
= - \boldsymbol d_i \times \big( \boldsymbol d_i\times\mathrm{Grad}_s\, \boldsymbol\omega \big)
= 2 \,\mathrm{Grad}_s\, \boldsymbol\omega\,,
\end{equation}
since 
\[  
- \boldsymbol d_i \times \big( \boldsymbol d_i\times \boldsymbol\omega,_\alpha \big) = - (\boldsymbol d_i\cdot \boldsymbol\omega,_\alpha ) \boldsymbol d_i + (\boldsymbol d_i\cdot\boldsymbol d_i) \boldsymbol\omega,_\alpha 
=-\boldsymbol\omega,_\alpha + 3\,\boldsymbol\omega,_\alpha
= 2\,\boldsymbol\omega,_\alpha\,.
\]
We substitute \eqref{f73} into \eqref{f72,5} and find
\begin{equation}\label{f74}
\delta\boldsymbol K^e = \boldsymbol Q_e^T\, \mathrm{Grad}_s \boldsymbol \omega.
\end{equation}
By virtue of \eqref{f69} and \eqref{f74}, the relation \eqref{f63} becomes
\begin{equation}\label{f75}
\delta\, 
W_{\mathrm{shell}} \, = \, \boldsymbol N  
: \big( \mathrm{Grad}_s (\delta \boldsymbol m) -  \boldsymbol\Omega \,\boldsymbol F_s \big) +
\boldsymbol M : \mathrm{Grad}_s \boldsymbol \omega\,.
\end{equation}
We can rewrite the term $ \boldsymbol N :(\boldsymbol \Omega\boldsymbol F_s) $ as follows
\begin{equation}\label{f76}
\boldsymbol N :(\boldsymbol \Omega\boldsymbol F_s) = - \boldsymbol \Omega :(\boldsymbol F_s\boldsymbol N^T) = -
\boldsymbol \omega\cdot \mathrm{axl} \big( \boldsymbol F_s \boldsymbol N^T- \boldsymbol N\boldsymbol F_s^T \big),
\end{equation}
since
\[  
\boldsymbol \Omega : \boldsymbol X = \mathrm{axl}(\boldsymbol \Omega) \cdot \mathrm{axl}(\boldsymbol X- \boldsymbol X^T)
\]
for any second order tensor $ \boldsymbol X $ and any skew-symmetric tensor $ \boldsymbol \Omega $. We use \eqref{f76} in \eqref{f75} and deduce
\begin{equation}\label{f77}
\delta\, 
W_{\mathrm{shell}} \, = \,  \boldsymbol N  
:  \mathrm{Grad}_s (\delta \boldsymbol m) +
\boldsymbol M : \mathrm{Grad}_s \boldsymbol \omega
+  \mathrm{axl} \big( \boldsymbol F_s \boldsymbol N^T- \boldsymbol N\boldsymbol F_s^T \big)\cdot \boldsymbol \omega\,.
\end{equation}
For the first two terms in the right-hand side of equation \eqref{f77} we employ relations of the type 
\[  
\boldsymbol S : \mathrm{Grad}_s \boldsymbol v = \mathrm{Div}_s(\boldsymbol S^T\boldsymbol v) - \big(\mathrm{Div}_s\boldsymbol S\big)\cdot \boldsymbol v\,,
\]
together with the divergence theorem on surfaces. Thus, in view of the null boundary conditions on $ \partial \omega_\xi $ we derive
\begin{equation}\label{f78}
\displaystyle\int_{\omega_\xi}  \boldsymbol N  :  \mathrm{Grad}_s (\delta \boldsymbol m)\, \mathrm{d}a =
\int_{\partial\omega_\xi}  (\delta \boldsymbol m)\cdot(\boldsymbol N  \boldsymbol \nu)\, \mathrm{d}\ell 
- \displaystyle\int_{\omega_\xi}  \big(\mathrm{Div}_s\boldsymbol N\big)  \cdot (\delta \boldsymbol m)\, \mathrm{d}a
= - \displaystyle\int_{\omega_\xi}  \big(\mathrm{Div}_s\boldsymbol N\big)  \cdot (\delta \boldsymbol m)\, \mathrm{d}a 
\end{equation}
and similarly
\begin{equation}\label{f79}
\displaystyle\int_{\omega_\xi}  \boldsymbol M  :  \mathrm{Grad}_s  \boldsymbol\omega\, \mathrm{d}a 
= - \displaystyle\int_{\omega_\xi}  \big(\mathrm{Div}_s\boldsymbol M\big)  \cdot  \boldsymbol \omega\, \mathrm{d}a .
\end{equation}
Finally, in view of \eqref{f77}-\eqref{f79} we obtain
\begin{equation}\label{f80}
0\,=\, \delta\, I \,=\,  \displaystyle\int_{\omega_\xi} \delta\, W_{\mathrm{shell}} \,\, \mathrm{d}a =
- \int_{\omega_\xi} \Big[ \big(\mathrm{Div}_s\boldsymbol N\big)  \cdot (\delta \boldsymbol m) + \big(\mathrm{Div}_s\boldsymbol M+ \mathrm{axl} \big(\boldsymbol N\boldsymbol F_s^T - \boldsymbol F_s \boldsymbol N^T\big)\big)  \cdot  \boldsymbol \omega \Big]\, \mathrm{d}a,
\end{equation}
for any virtual translation $ \delta\boldsymbol m $ and any virtual rotation $ \boldsymbol\omega = \mathrm{axl} \big((\delta\boldsymbol Q_e)\boldsymbol Q_e^T\big) $. Relation \eqref{f80} yields the following local forms of the equilibrium equations
\begin{equation}\label{f81}
\mathrm{Div}_s\boldsymbol N  = \boldsymbol 0\qquad\mathrm{and}\qquad 
\mathrm{Div}_s\boldsymbol M + \mathrm{axl} \big(\boldsymbol N\boldsymbol F_s^T - \boldsymbol F_s \boldsymbol N^T\big)= \boldsymbol 0.
\end{equation}

\textbf{Remark:} The principle of virtual work for 6-parameter shells corresponding to equation 
\eqref{f80} has been presented in \cite{Eremeyev06,Birsan-Neff-L58-2017}. \hfill$ \Box $\medskip

If we consider now external body forces $ \boldsymbol f $ and couples $ \boldsymbol c $\,, we can write the equilibrium equations for Cosserat shells in the general form (see, e.g. \cite{Eremeyev06,Birsan-Neff-L58-2017})
\begin{equation}\label{f82}
\mathrm{Div}_s\boldsymbol N + \boldsymbol f = \boldsymbol 0,\qquad 
\mathrm{Div}_s\boldsymbol M + \mathrm{axl} \big(\boldsymbol N\boldsymbol F_s^T - \boldsymbol F_s \boldsymbol N^T\big) + \boldsymbol c= \boldsymbol 0. 
\end{equation}
The tensors $ \boldsymbol N $ and $ \boldsymbol M $ are the internal surface stress tensor and the internal surface couple tensor (of the first Piola-Kirchhoff type), respectively. They are given by the relations \eqref{f64}.

The general form of the boundary conditions of mixed type on $ \partial \omega_\xi $ is (see, e.g. \cite{Pietraszkiewicz04,Pietraszkiewicz11,Birsan-Neff-MMS-2014})
\begin{equation}\label{f83}
\begin{array}{rcl}
\boldsymbol{N}\boldsymbol{\nu} & = & \boldsymbol{N}^*,\qquad \boldsymbol{M}\boldsymbol{\nu}  =  \boldsymbol{M}^*\quad\mathrm{along}\,\,\,\partial \omega_f\,,
\vspace{4pt}\\
\,\,\boldsymbol{m} & = & \boldsymbol{m}^* ,\qquad\quad \boldsymbol{Q}_e  =  \boldsymbol{Q}^* \quad\mathrm{along}\,\,\,\partial \omega_d\,,
\end{array}
\end{equation}
where $\partial \omega_f$ and $ \partial \omega_d $ build a disjoint partition  of the boundary curve $\partial \omega_\xi\,$. Here, $\boldsymbol{N}^*$ and $\boldsymbol{M}^*$ are the external boundary force and couple vectors respectively, applied along the deformed boundary curve, but measured per unit length of $\partial \omega_f\,$. On the portion of the boundary $\partial \omega_d$ we have Dirichlet-type boundary conditions for the deformation vector $ \boldsymbol{m} $ and the microrotation tensor $ \boldsymbol{Q}_e\, $.
\smallskip

Using the obtained form of the energy density \eqref{f61} and the relations \eqref{f64}, we can give the stress-strain relations in explicit form for our shell model. These will be written in the next section.

\section{Remarks and discussions on the Cosserat shell model} \label{Sect5}

In this section we write the strain-energy density \eqref{f61} in some alternative useful  forms and give the explicit expression for the constitutive equations \eqref{f64}. This allows us to compare the derived shell model with other approaches to 6-parameter shells and with the classical Koiter shell model.

We notice that the shell strain measures $ \boldsymbol E^e $ and $ \boldsymbol K^e $ (as well as the shell stress tensors $ \boldsymbol{Q}_e^T\boldsymbol N $ and $ \boldsymbol{Q}_e^T\boldsymbol M $) are tensors of the form $ \boldsymbol X=X_{i\alpha}\boldsymbol a^i\otimes \boldsymbol a^\alpha $ (where $ \boldsymbol a^3=\boldsymbol n_0 $). In what follows, we shall decompose any such tensor $ \boldsymbol X=X_{i\alpha}\boldsymbol a^i\otimes \boldsymbol a^\alpha $ in its ``planar'' part $ \boldsymbol a \boldsymbol X = X_{\beta\alpha}\boldsymbol a^\beta\otimes \boldsymbol a^\alpha  $ and its ``transversal'' part $ \boldsymbol n_0\boldsymbol X = X_{3\alpha}\boldsymbol a^\alpha$ according to
\begin{equation}\label{f84}
\boldsymbol X = \id_3 \boldsymbol X = (\boldsymbol a + \boldsymbol n_0 \otimes \boldsymbol n_0)\boldsymbol X =\boldsymbol a \boldsymbol X + \boldsymbol n_0 \otimes (\boldsymbol n_0 \boldsymbol X).
\end{equation}
Note that $ \boldsymbol a \boldsymbol X $ is a planar tensor in the tangent plane, while $ \boldsymbol n_0 \boldsymbol X $ is a vector in the tangent plane. For instance, the decomposition of the shell strain tensor $ \boldsymbol E^e  $ yields
\begin{equation}\label{f85}
\boldsymbol E^e= 
\boldsymbol a \boldsymbol E^e + \boldsymbol n_0 \otimes (\boldsymbol n_0 \boldsymbol E^e), \qquad \boldsymbol a\boldsymbol E^e= E^e_{\beta\alpha}\boldsymbol a^\beta\otimes \boldsymbol a^\alpha,\qquad 
\boldsymbol n_0 \boldsymbol E^e = E^e_{3\alpha}\boldsymbol a^\alpha,
\end{equation}
where $ \boldsymbol n_0 \boldsymbol E^e  $ describes the transverse shear deformations and $ \boldsymbol a\boldsymbol E^e $ the in-plane deformation of the shell. 

With this representation, we can decompose the constitutive equations 
\eqref{f64} in the following way
\begin{equation}\label{f87}
\boldsymbol a \boldsymbol Q_e^T  \boldsymbol N =  \dfrac{\partial W_{\mathrm{shell}}}{\partial (\boldsymbol a \boldsymbol E^e)}\;,\quad
\boldsymbol n_0 \boldsymbol Q_e^T  \boldsymbol N =  \dfrac{\partial W_{\mathrm{shell}}}{\partial (\boldsymbol n_0 \boldsymbol E^e)}\;,\quad
\boldsymbol a \boldsymbol Q_e^T  \boldsymbol M =  \dfrac{\partial W_{\mathrm{shell}}}{\partial (\boldsymbol a \boldsymbol K^e)}\;,\quad
\boldsymbol n_0 \boldsymbol Q_e^T  \boldsymbol M =  \dfrac{\partial W_{\mathrm{shell}}}{\partial (\boldsymbol n_0 \boldsymbol K^e)}\;.
\end{equation}

\subsection{Explicit stress-strain relations} \label{Sect5.1}

In order to write the stress-strain relations explicitly, let us put the equations \eqref{f50} and \eqref{f51} in the forms
\begin{equation}\label{f88}
\begin{array}{rl}
W_{\mathrm{Coss}}(\boldsymbol X,\boldsymbol Y) =& 
\mu\, \mathrm{sym}( \boldsymbol a\boldsymbol X) : \mathrm{sym} (\boldsymbol a \boldsymbol Y) +  \mu_c\, \mathrm{skew}(\boldsymbol a \boldsymbol X) : \mathrm{skew}(\boldsymbol a \boldsymbol Y) +\,\dfrac{\lambda\,\mu}{\lambda+2\mu}\,\big( \mathrm{tr} \boldsymbol X\big)\,\big(\mathrm{tr} \boldsymbol Y\big) 
\vspace{6pt}\\
&+ \,\dfrac{2\mu\,\mu_c}{\mu+\mu_c}\,\big(  \boldsymbol n_0\boldsymbol X\big)\cdot\big(\boldsymbol n_0\boldsymbol Y\big),
\vspace{6pt}\\
W_{\mathrm{Coss}}(\boldsymbol X)  = &
\mu\,\| \mathrm{sym}( \boldsymbol a\boldsymbol X) \|^2 +  \mu_c\| \mathrm{skew}(\boldsymbol a \boldsymbol X) \|^2 +\,\dfrac{\lambda\,\mu}{\lambda+2\mu}\,\big( \mathrm{tr} \boldsymbol X\big)^2 
+ \,\dfrac{2\mu\,\mu_c}{\mu+\mu_c}\,\| \boldsymbol n_0\boldsymbol X\|^2
\end{array}
\end{equation}
and note that $ \mathrm{tr} \boldsymbol X= \mathrm{tr} (\boldsymbol a\boldsymbol X) $. Suggested by \eqref{f88}, we introduce the fourth order planar tensor  $ \underline{\boldsymbol C}_S $ of elastic moduli for the shell
\begin{equation}\label{f89}
\begin{array}{c}
\underline{\boldsymbol C}_S = C_S^{\alpha\beta\gamma\delta} \boldsymbol a_\alpha\otimes\boldsymbol a_\beta\otimes\boldsymbol a_\gamma\otimes\boldsymbol a_\delta\qquad \mathrm{with}
\vspace{6pt}\\
 C_S^{\alpha\beta\gamma\delta} = \mu\,\big( a^{\alpha\gamma}a^{\beta\delta} + a^{\alpha\delta}a^{\beta\gamma}  \big)
+ \mu_c\,\big( a^{\alpha\gamma}a^{\beta\delta} - a^{\alpha\delta}a^{\beta\gamma}  \big)
+ \,\dfrac{2\lambda\,\mu}{\lambda+2\mu}\; a^{\alpha\beta}a^{\gamma\delta} \,.
\end{array}
\end{equation}
Then, the tensor $ \underline{\boldsymbol C}_S $ satisfies the major symmetries $  C_S^{\alpha\beta\gamma\delta} =  C_S^{\gamma\delta\alpha\beta}  $ and we have
\begin{equation}\label{f90}
\underline{\boldsymbol C}_S : \boldsymbol T = 2\mu\, \mathrm{sym}\, \boldsymbol T  + 2 \mu_c\, \mathrm{skew}\,\boldsymbol T +\,\dfrac{2\lambda\,\mu}{\lambda+2\mu}\,\big( \mathrm{tr}\, \boldsymbol T\big)\,\boldsymbol a\,,
\end{equation}
for any planar tensor $ \boldsymbol T = T_{\alpha\beta} \boldsymbol a^\alpha\otimes \boldsymbol a^\beta$. Due to the symmetry, the relations 
\eqref{f88} can be written in a simple way
\begin{equation}\label{f91}
\begin{array}{rl}
W_{\mathrm{Coss}}(\boldsymbol X,\boldsymbol Y) =& 
\dfrac12( \boldsymbol a\boldsymbol X) : \underline{\boldsymbol C}_S : ( \boldsymbol a\boldsymbol Y)+ \,\dfrac{2\mu\,\mu_c}{\mu+\mu_c}\,\big(  \boldsymbol n_0\boldsymbol X\big)\cdot\big(\boldsymbol n_0\boldsymbol Y\big)
=\, \dfrac12 \, C_S^{\alpha\beta\gamma\delta} X_{\alpha\beta}Y_{\gamma\delta} + \dfrac{2\mu\,\mu_c}{\mu+\mu_c}\,X_{3\alpha}Y_{3\alpha},
\vspace{6pt}\\
W_{\mathrm{Coss}}(\boldsymbol X)  = &
\dfrac12( \boldsymbol a\boldsymbol X) : \underline{\boldsymbol C}_S : ( \boldsymbol a\boldsymbol X)+ \,\dfrac{2\mu\,\mu_c}{\mu+\mu_c}\,\| \boldsymbol n_0\boldsymbol X\|^2,
\end{array}
\end{equation}
for any tensors $ \boldsymbol X=X_{i\alpha}\boldsymbol a^i\otimes \boldsymbol a^\alpha $, $ \boldsymbol Y=Y_{i\alpha}\boldsymbol a^i\otimes \boldsymbol a^\alpha $.

Similarly, the quadratic form $ W_{\mathrm{curv}} $ defined by \eqref{f7} can be put in the form
\begin{equation}\label{f92}
\begin{array}{rl}
W_{\mathrm{curv}}(\boldsymbol X) =& \mu\,L_c^2\Big(\,
b_1\| \mathrm{sym}( \boldsymbol a\boldsymbol X) \|^2 +  b_2\| \mathrm{skew}(\boldsymbol a \boldsymbol X) \|^2 +\,\big(b_3-\,\dfrac{b_1}{3}\big)\big( \mathrm{tr} \boldsymbol X\big)^2 
+ \,\dfrac{b_1+b_2}{2}\,\| \boldsymbol n_0\boldsymbol X\|^2
\Big)
\vspace{6pt}\\
  = &
\dfrac12( \boldsymbol a\boldsymbol X) : \underline{\boldsymbol G}_S : ( \boldsymbol a\boldsymbol X)+ \mu\,L_c^2\,\dfrac{b_1+b_2}{2}\,\| \boldsymbol n_0\boldsymbol X\|^2 
\,= \dfrac12 \, G_S^{\alpha\beta\gamma\delta} X_{\alpha\beta}X_{\gamma\delta} +  \mu\,L_c^2\,\dfrac{b_1+b_2}{2}\,X_{3\alpha}X_{3\alpha}
\end{array}
\end{equation}
for any tensor $ \boldsymbol X=X_{i\alpha}\boldsymbol a^i\otimes \boldsymbol a^\alpha $, where the fourth order planar tensor $ \underline{\boldsymbol G}_S $ is given by
\begin{equation}\label{f93}
\begin{array}{c}
\underline{\boldsymbol G}_S = G_S^{\alpha\beta\gamma\delta} \boldsymbol a_\alpha\otimes\boldsymbol a_\beta\otimes\boldsymbol a_\gamma\otimes\boldsymbol a_\delta\qquad \mathrm{with}
\vspace{6pt}\\
G_S^{\alpha\beta\gamma\delta} = \mu\,L_c^2\Big( 
b_1\,\big( a^{\alpha\gamma}a^{\beta\delta} + a^{\alpha\delta}a^{\beta\gamma}  \big)
+ b_2\,\big( a^{\alpha\gamma}a^{\beta\delta} - a^{\alpha\delta}a^{\beta\gamma}  \big)
+ \big(b_3-\,\dfrac{b_1}{3}\,\big)\, a^{\alpha\beta}a^{\gamma\delta} \,
\Big).
\end{array}
\end{equation}
We see that $ G_S^{\alpha\beta\gamma\delta} =  G_S^{\gamma\delta\alpha\beta} $ and for any planar tensor $ \boldsymbol T = T_{\alpha\beta} \boldsymbol a^\alpha\otimes \boldsymbol a^\beta$ it holds
\begin{equation}\label{f94}
\underline{\boldsymbol G}_S : \boldsymbol T = 2\mu\,L_c^2 \Big(b_1\,\mathrm{sym}\, \boldsymbol T  + b_2\, \mathrm{skew}\,\boldsymbol T +\big(b_3-\,\dfrac{b_1}{3}\big)\big( \mathrm{tr}\, \boldsymbol T\big)\,\boldsymbol a\Big).
\end{equation}
In order to show that the quadratic forms $ W_{\mathrm{Coss}} $ and $ W_{\mathrm{curv}} $ are positive definite, let us introduce the \textit{surface deviator operator} $ \mathrm{dev}_s $ defined by \cite{Birsan-Neff-L58-2017}
\begin{equation}\label{f95}
\mathrm{dev}_s\, \boldsymbol X \,:=\, \boldsymbol X \,-\,
\dfrac12 \big( \mathrm{tr}\, \boldsymbol X\big)\,\boldsymbol a.
\end{equation}
According to Lemma 2.1 in \cite{Birsan-Neff-L58-2017} we can decompose any tensor $ \boldsymbol X=X_{i\alpha}\boldsymbol a^i\otimes \boldsymbol a^\alpha $ as a \textit{direct sum} (orthogonal decomposition) as follows
\begin{equation}\label{f96}
\boldsymbol{X}\,=\, \mathrm{dev_ssym}\, \boldsymbol{X}\, + \, \mathrm{skew}\, \boldsymbol{X}\,   + \, \frac{1}{2}\,\big(\mathrm{tr}\,\boldsymbol{X}\big)\,\boldsymbol a\,.
\end{equation}
Then, relations \eqref{f95} and \eqref{f96} imply 
\begin{equation}\label{f97}
\mathrm{sym}\,\boldsymbol{X}\,=\, \mathrm{dev_ssym}\, \boldsymbol{X}\,    + \, \frac{1}{2}\,\big(\mathrm{tr}\,\boldsymbol{X}\big)\,\boldsymbol a
\qquad \mathrm{and}\qquad 
\|\mathrm{sym}\,\boldsymbol{X}\|^2=\| \mathrm{dev_ssym}\, \boldsymbol{X}\|^2    + \, \frac{1}{2}\,\big(\mathrm{tr}\,\boldsymbol{X}\big)^2.
\end{equation}
Substituting \eqref{f97} into the relations \eqref{f90} and \eqref{f94}, we get (for any $ \boldsymbol T = T_{\alpha\beta} \boldsymbol a^\alpha\otimes \boldsymbol a^\beta$)
\begin{equation}\label{f98}
\begin{array}{l}
\underline{\boldsymbol C}_S : \boldsymbol T = 2\mu\, \mathrm{dev_ssym}\, \boldsymbol T  + 2 \mu_c\, \mathrm{skew}\,\boldsymbol T +\,\dfrac{\mu(3\lambda+2\mu)}{\lambda+2\mu}\,\big( \mathrm{tr}\, \boldsymbol T\big)\,\boldsymbol a\,,
\vspace{6pt}\\
\underline{\boldsymbol G}_S : \boldsymbol T = 2\mu\,L_c^2 \Big(b_1\,\mathrm{dev_ssym}\, \boldsymbol T  + b_2\, \mathrm{skew}\,\boldsymbol T +\big(b_3+\,\dfrac{b_1}{6}\big)\big( \mathrm{tr}\, \boldsymbol T\big)\,\boldsymbol a\Big)
\end{array}
\end{equation}
and the quadratic forms \eqref{f88}$ _2 $ and \eqref{f92} become
\begin{equation}\label{f99}
\begin{array}{rl}
W_{\mathrm{Coss}}(\boldsymbol X)  = &
\mu\,\| \mathrm{dev_ssym}( \boldsymbol a\boldsymbol X) \|^2 +  \mu_c\| \mathrm{skew}(\boldsymbol a \boldsymbol X) \|^2 +\,\dfrac{\mu(3\lambda+2\mu)}{2(\lambda+2\mu)}\,\big( \mathrm{tr} \boldsymbol X\big)^2 
+ \,\dfrac{2\mu\,\mu_c}{\mu+\mu_c}\,\| \boldsymbol n_0\boldsymbol X\|^2 ,
\vspace{6pt}\\
W_{\mathrm{curv}}(\boldsymbol X) =& \mu\,L_c^2\Big(\,
b_1\| \mathrm{dev_ssym}( \boldsymbol a\boldsymbol X) \|^2 +  b_2\| \mathrm{skew}(\boldsymbol a \boldsymbol X) \|^2 +\,\big(b_3+\,\dfrac{b_1}{6}\big)\big( \mathrm{tr} \boldsymbol X\big)^2 
+ \,\dfrac{b_1+b_2}{2}\,\| \boldsymbol n_0\boldsymbol X\|^2
\Big).
\end{array}
\end{equation}
Under the usual assumptions on the material constants $ \mu>0 $, $ 3\lambda+2\mu>0 $ (from classical elasticity), together with $ \mu_c>0 $ and $ b_i>0 $, we see now that the quadratic forms \eqref{f99} are positive definite, since all the coefficients are positive.

Finally, we substitute \eqref{f91}, \eqref{f92} in the strain-energy density \eqref{f61} and performing the differentiation according to the relations \eqref{f87}, we obtain the following explicit forms of the constitutive equations for
the internal surface stress tensor $\boldsymbol Q_e^T  \boldsymbol N $ and the internal surface couple tensor  $ \boldsymbol Q_e^T \boldsymbol M $ of Cosserat shells
\[  
\boldsymbol Q_e^T  \boldsymbol N = 
\boldsymbol a \boldsymbol Q_e^T  \boldsymbol N + \boldsymbol n_0 \otimes (\boldsymbol n_0 \boldsymbol Q_e^T  \boldsymbol N),\qquad
\boldsymbol Q_e^T  \boldsymbol M = 
\boldsymbol a \boldsymbol Q_e^T  \boldsymbol M + \boldsymbol n_0 \otimes (\boldsymbol n_0 \boldsymbol Q_e^T  \boldsymbol M)
\]
with
\begin{equation}\label{f100}
\begin{array}{rl}
\boldsymbol a \boldsymbol Q_e^T  \boldsymbol N = &  
\Big(h- K\, \dfrac{h^3}{12}\Big) \,
\underline{\boldsymbol C}_S : \big(\boldsymbol a \boldsymbol E^e\big)
+ \,\dfrac{h^3}{12}\,\big[\,
\underline{\boldsymbol C}_S : \big(\boldsymbol a \boldsymbol E^e\boldsymbol b + \boldsymbol c \boldsymbol K^e\big)
\big]\boldsymbol b
-\dfrac{h^3}{12}\; \underline{\boldsymbol C}_S : \big( \boldsymbol c \boldsymbol K^e\boldsymbol b^*\big),
\vspace{6pt}\\
\boldsymbol n_0 \boldsymbol Q_e^T  \boldsymbol N =  & 
\,\dfrac{4\mu\,\mu_c}{\mu+\mu_c}\,\Big[
\Big(h- 2K\,\dfrac{h^3}{12}\,\Big) \big(\boldsymbol n_0 \boldsymbol E^e\big) + 2H\,\dfrac{h^3}{12}\, \big(\boldsymbol n_0 \boldsymbol E^e\boldsymbol b\big)
\Big]
 ,
\vspace{6pt}\\
\boldsymbol a \boldsymbol Q_e^T  \boldsymbol M = & \Big(h- K\, \dfrac{h^3}{12}\Big) \,
\underline{\boldsymbol G}_S : \big(\boldsymbol a \boldsymbol K^e\big)
+ \,\dfrac{h^3}{12}\,\boldsymbol c\,\big[\,
\underline{\boldsymbol C}_S : \big(\boldsymbol a \boldsymbol E^e\boldsymbol b + \boldsymbol c \boldsymbol K^e\big)
\big]
-\dfrac{h^3}{12}\,\boldsymbol c\, \big[\underline{\boldsymbol C}_S : \big( \boldsymbol a \boldsymbol E^e\big)\big]\boldsymbol b^* 
\vspace{6pt}\\
& + \dfrac{h^3}{12}\,\, \big[\underline{\boldsymbol G}_S : \big( \boldsymbol a \boldsymbol K^e\boldsymbol b\big)\big]\boldsymbol b\,,
\vspace{6pt}\\
\boldsymbol n_0 \boldsymbol Q_e^T  \boldsymbol M =  & 
\mu\,L_c^2 \,(b_1+b_2) \Big[
\Big(h- 2K\,\dfrac{h^3}{12}\,\Big) \big(\boldsymbol n_0 \boldsymbol K^e\big) + 2H\,\dfrac{h^3}{12} \, \big(\boldsymbol n_0 \boldsymbol K^e\boldsymbol b\big)
\Big]
,
\end{array}
\end{equation}
where the tensors of elastic moduli $ \underline{\boldsymbol C}_S $ and $ \underline{\boldsymbol G}_S $ are given in \eqref{f89}, \eqref{f90} and \eqref{f93}, \eqref{f94}.

\subsection{Comparison with other 6-parameter shell models}\label{Sect5.2}

We present a detailed comparison with the 
related shell model of order $ O(h^5) $ which has been presented recently in \cite{Birsan-Neff-MMS-2019}. The Cosserat shell model derived in \cite{Birsan-Neff-MMS-2019} has many similarities with the present model, but there are also some differences, which we indicate now.

First of all, the derivation method and starting point in \cite{Birsan-Neff-MMS-2019} is different, since the deformation function $ \boldsymbol \varphi $ is assumed to be quadratic in $ x_3\, $. More precisely, the following ansatz is adopted (see \cite[f. (65)]{Birsan-Neff-MMS-2019})
\begin{equation}\label{f101}
\boldsymbol \varphi( x_i) = \boldsymbol m(x_1,x_2)+ x_3\,\alpha(x_1,x_2)\, \boldsymbol d_3 +\dd\frac{x_3^2}{2}\,\beta(x_1,x_2) \,\boldsymbol d_3\, .
\end{equation}
If we compare this ansatz with the expansion \eqref{f27}, we see the assumption  \eqref{f101} is more restrictive.

On the other hand, the hypotheses \eqref{f42} from the classical shell theory were replaced in \cite{Birsan-Neff-MMS-2019}  by the weaker requirements (see \cite[f. (60)]{Birsan-Neff-MMS-2019})
\begin{equation}\label{f102}
\boldsymbol n_0\cdot \boldsymbol T_0 \boldsymbol n_0 =   0 \qquad\mathrm{and}\qquad 
\boldsymbol n_0\cdot \boldsymbol T_0^{\,\prime} \boldsymbol n_0 =  0 ,
\end{equation}
i.e. only the normal components of the stress vectors $ \boldsymbol t^+ $ , $ \boldsymbol t^- $ on the upper and lower surfaces of the shell are assumed to be zero. The two scalar equations \eqref{f102} are then employed in \cite{Birsan-Neff-MMS-2019} to determine the two scalar coefficients $ \alpha(x_1,x_2) $ and $ \beta(x_1,x_2) $ appearing in \eqref{f101}.
Moreover, we note that the paper \cite{Birsan-Neff-MMS-2019} presents a shell model of order $ O(h^5) $.

This different approach leads to a slightly different form of the strain-energy density. If we retain only the terms up to the  order $ O(h^3) $ in the strain-energy density (see \cite[f. (104)]{Birsan-Neff-MMS-2019}), we get
\begin{equation}\label{f103}
\begin{array}{rl}
\widehat{W}_{\mathrm{shell}}(\boldsymbol E^e,\boldsymbol K^e) =  &
\Big(h- K\, \dfrac{h^3}{12}\Big) \big[  W_{\mathrm{mixt}}\big( \boldsymbol E^e\big) + W_{\mathrm{curv}} \big(\boldsymbol K^e\big) \big]
\vspace{6pt}\\
& + \;\dfrac{h^3}{12}\,\big[  W_{\mathrm{mixt}}\big(\boldsymbol E^e\boldsymbol b + \boldsymbol c \boldsymbol K^e \big) 
-2 W_{\mathrm{mixt}}\big(\boldsymbol E^e\,,\, \boldsymbol c \boldsymbol K^e \boldsymbol b^*  \big)
+ W_{\mathrm{curv}} \big(\boldsymbol K^e \boldsymbol b\big) \big],
\end{array}
\end{equation}
where $ {W}_{\mathrm{mixt}} $ is given by \eqref{f49,5}. We compare this expression with our energy \eqref{f61}.
Using the decomposition of tensors in planar and transversal parts \eqref{f84}, we deduce from 
\eqref{f49,5} and \eqref{f50} the relations
\begin{equation}\label{f104}
\begin{array}{l}
W_{\mathrm{mixt}}(\boldsymbol S,\boldsymbol T) =   W_{\mathrm{mixt}}(\boldsymbol a\boldsymbol S,\boldsymbol a\boldsymbol T)+ \, \dfrac{\mu+\mu_c}{2}\,\big(\boldsymbol n_0  \boldsymbol S\big)\cdot\big(\boldsymbol n_0  \boldsymbol T\big),
	\vspace{6pt}\\
W_{\mathrm{Coss}}(\boldsymbol S,\boldsymbol T) =  W_{\mathrm{mixt}}(\boldsymbol a\boldsymbol S,\boldsymbol a\boldsymbol T)+ \, \dfrac{2\mu\,\mu_c}{\mu+\mu_c}\,\big(\boldsymbol n_0  \boldsymbol S\big)\cdot\big(\boldsymbol n_0  \boldsymbol T\big).
\end{array}
\end{equation}
Thus, using the relation  \eqref{f104}$ _1 $ the strain-energy density \eqref{f103} (obtained in \cite{Birsan-Neff-MMS-2019} for order $ O(h^3) $) becomes
\begin{equation}\label{f106}
\begin{array}{l}
\widehat{W}_{\mathrm{shell}}(\boldsymbol E^e,\boldsymbol K^e) =   
\Big(h- K\, \dfrac{h^3}{12}\Big) \big[  W_{\mathrm{mixt}}\big( \boldsymbol a \boldsymbol E^e\big) 
+ \, \dfrac{\mu+\mu_c}{2}\,\|\boldsymbol n_0  \boldsymbol E^e \|^2
+ W_{\mathrm{curv}} \big(\boldsymbol K^e\big) \big]
\vspace{6pt}\\
\qquad\qquad + \;\dfrac{h^3}{12}\,\big[  W_{\mathrm{mixt}}\big(\boldsymbol a \boldsymbol E^e\boldsymbol b 
+ \boldsymbol c \boldsymbol K^e \big)
+ \, \dfrac{\mu+\mu_c}{2}\,\|\boldsymbol n_0  \boldsymbol E^e\boldsymbol b \|^2 
-2 W_{\mathrm{mixt}}\big(\boldsymbol a \boldsymbol E^e\,,\, \boldsymbol c \boldsymbol K^e \boldsymbol b^*  \big)
+ W_{\mathrm{curv}} \big(\boldsymbol K^e \boldsymbol b\big) \big].
\end{array}
\end{equation}
On the other hand, our strain-energy density \eqref{f61} can be written with the help of \eqref{f104}$ _2 $ in the following alternative form
\begin{equation}\label{f107}
\begin{array}{l}
W_{\mathrm{shell}}(\boldsymbol E^e,\boldsymbol K^e) =  
\Big(h- K\, \dfrac{h^3}{12}\Big) \big[  W_{\mathrm{mixt}}\big( \boldsymbol a \boldsymbol E^e\big) 
+\, \dfrac{2\mu\,\mu_c}{\mu+\mu_c}\,\|\boldsymbol n_0  \boldsymbol E^e \|^2
+ W_{\mathrm{curv}} \big(\boldsymbol K^e\big) \big]
\vspace{6pt}\\
\qquad\qquad + \;\dfrac{h^3}{12}\,\big[  W_{\mathrm{mixt}}\big(\boldsymbol a \boldsymbol E^e\boldsymbol b 
+ \boldsymbol c \boldsymbol K^e \big)
+ \, \dfrac{2\mu\,\mu_c}{\mu+\mu_c}\,\|\boldsymbol n_0  \boldsymbol E^e\boldsymbol b \|^2 
-2 W_{\mathrm{mixt}}\big(\boldsymbol a \boldsymbol E^e\,,\, \boldsymbol c \boldsymbol K^e \boldsymbol b^*  \big)
+ W_{\mathrm{curv}} \big(\boldsymbol K^e \boldsymbol b\big) \big].
\end{array}
\end{equation}
By comparison of \eqref{f106} and \eqref{f107} we see that the only difference between these two 
strain-energy densities resides in the coefficients of the transverse shear deformation terms $ \| \boldsymbol n_0  \boldsymbol E^e \|^2$ and $ \| \boldsymbol n_0  \boldsymbol E^e \boldsymbol b \|^2 $. All other terms and coefficients in \eqref{f106} and \eqref{f107} are identical.

Note that the transverse shear coefficient in the present model \eqref{f107} is the harmonic mean $ \, \dfrac{2\mu\,\mu_c}{\mu+\mu_c}\, $ , while in the energy density \eqref{f106} (derived in \cite{Birsan-Neff-MMS-2019})  it is the arithmetic mean $ \, \dfrac{\mu+\mu_c}{2}\, $ . We mention that the same coefficient $ \, \dfrac{2\mu\,\mu_c}{\mu+\mu_c}\, $ for the transverse shear energy has been obtained using $ \Gamma $-convergence in \cite{Neff_Hong_Reissner08} in the case of plates. 
This confirms the result \eqref{f107} obtained in our present work. We remind that  this coefficient is adjusted in many plate and shell models by a correction factor, the so-called \textit{shear correction factor} (see for instance the discussions in  \cite{Altenbach00,Pietraszkiewicz10,Vlachoutsis}).
\medskip

\textbf{Further remarks:}

1. We remark that the strain-energy density \eqref{f107} obtained in this paper satisfies the invariance properties required by the local symmetry group of isotropic 6-parameter shells. These invariance requirements have been established in a general theoretical framework in  \cite[Section 9]{Eremeyev06}. 

2. The form of the constitutive relation \eqref{f107} (equivalent to \eqref{f61}) is remarkable, since one cannot find in the literature on 6-parameter shells appropriate expressions of the 
strain-energy density $ W_{\mathrm{shell}}(\boldsymbol E^e,\boldsymbol K^e) $ with coefficients depending on the initial curvature $ \boldsymbol b $ and expressed in terms of the three-dimensional material constants. Indeed, the strain-energy densities proposed in the literature are  either simple expressions with constant coefficients (see, e.g. \cite[f. (72)]{Birsan-Neff-MMS-2014}, \cite[f. (50)]{Birsan-Neff-L58-2017}, \cite{Pietraszkiewicz-book04,Pietraszkiewicz10}), or general quadratic forms of $ \boldsymbol E^e,\;\boldsymbol K^e $ with unidentified coefficients (see, e.g. \cite[f. (52)]{Eremeyev06}.

3. We mention that the numerical treatment for the related \textit{planar} Cosserat shell model derived in \cite{Neff_plate04_cmt,Neff_plate07_m3as} has been presented in \cite{Sander-Neff-Birsan-16}, using geodesic finite elements.

4. If the thickness $ h $ is sufficiently small, one can show that the strain-energy density $ W_{\mathrm{shell}}(\boldsymbol E^e,\boldsymbol K^e) $ is a coercive and convex function of its arguments. Then, in view of Theorem 6 from \cite{Birsan-Neff-MMS-2014}, one can prove the existence of minimizers for our nonlinear Cosserat shell model.

\subsection{Relation to the classical Koiter shell model}\label{Sect5.3}

In this section, we discuss the relation to the classical shell theory and show that our strain-energy density \eqref{f107} can be reduced, in a certain sense, to the strain-energy of the classical Koiter model.

Thus, if we consider that the three-dimensional material is a Cauchy continuum (with no microrotation), then the Cosserat couple modulus and the curvature energy $ W_{\mathrm{curv}} $ are vanishing in the model 
\eqref{f5}-\eqref{f6}:
\begin{equation}\label{f108}
\mu_c =   0, \qquad W_{\mathrm{curv}}\equiv  0 .
\end{equation}
Hence, the fourth order constitutive tensor for shells \eqref{f89} reduces to
\begin{equation}\label{f109}
   C_S^{\alpha\beta\gamma\delta} = \mu\,\big( a^{\alpha\gamma}a^{\beta\delta} + a^{\alpha\delta}a^{\beta\gamma}  \big)
 + \,\dfrac{2\lambda\,\mu}{\lambda+2\mu}\; a^{\alpha\beta}a^{\gamma\delta}\,,
\end{equation}
which coincide with the tensor of linear plane-stress elastic moduli, that appears in the Koiter model (see, e.g. \cite{Koiter60}, \cite[Sect. 4.1]{Ciarlet05}, \cite[f. (101)]{Steigmann13}). In view of \eqref{f49,6} and \eqref{f108}$ _1\, $, we notice that in this case
\begin{equation}\label{f110}
W_{\mathrm{mixt}}(\boldsymbol S)= W_{\mathrm{Koit}}( \boldsymbol S),
\end{equation}
where
\begin{equation}\label{f110,5}
W_{\mathrm{Koit}}(\boldsymbol S) := \mu\,\|\mathrm{sym}\, \boldsymbol S\|^2 + \dfrac{\lambda\, \mu}{\lambda+2 \mu}\,(\mathrm{tr}\, \boldsymbol S)^2
\end{equation}
is the quadratic form appearing in the Koiter model. We remind that the areal strain-energy density for Koiter shells has the expression \cite{Koiter60,Ciarlet05,Steigmann13}
\begin{equation}\label{f111}
 h\,W_{\mathrm{Koit}}(\boldsymbol \epsilon) + \dfrac{h^3}{12}\,W_{\mathrm{Koit}}(\boldsymbol \rho),
\end{equation}
where the \textit{change of metric} tensor $ \boldsymbol \epsilon $ and the \textit{change of curvature} tensor $ \boldsymbol \rho $ are the nonlinear shell strain measures, which are given by
\begin{equation}\label{f112}
\begin{array}{l}
\boldsymbol \epsilon \,=\,  \dfrac12\,\big( \boldsymbol m,_\alpha \cdot\, \boldsymbol m,_\beta -a_{\alpha\beta}\big) \,\boldsymbol a^\alpha \otimes \boldsymbol a^\beta\, =\, \dfrac12\,\big[\, (\mathrm{Grad}_s\boldsymbol m)^T  (\mathrm{Grad}_s\boldsymbol m) - \boldsymbol a\,\big],
\vspace{6pt}\\
\boldsymbol \rho \,= \,\big( \boldsymbol n\cdot \boldsymbol m,_{\alpha\beta} - \boldsymbol n_0\cdot \boldsymbol a_{\alpha,\beta}\big) \,\boldsymbol a^\alpha \otimes \boldsymbol a^\beta\, =\, - (\mathrm{Grad}_s\boldsymbol m)^T  (\mathrm{Grad}_s\boldsymbol n) - \boldsymbol b\,.
\end{array}
\end{equation}
Here, $ \boldsymbol n $ designates the unit normal vector to the deformed midsurface and we note that $ \boldsymbol \epsilon $ and  $ \boldsymbol \rho $ are symmetric planar tensors.

To obtain the classical shell model as a special case of our approach, we adopt the Kirchhoff-Love hypotheses. Thus, we assume that the reference unit normal $ \boldsymbol n_0 $ becomes after deformation the unit normal to the deformed midsurface, i.e. $ \boldsymbol n_0 $ transforms to $ \boldsymbol n $. But since we have $  \boldsymbol Q_e\boldsymbol n_0  = \boldsymbol Q_e\boldsymbol d_3^0 = \boldsymbol d_3\, $, this assumption means that
\begin{equation}\label{f113}
\boldsymbol n = \boldsymbol d_3\,.
\end{equation}
Then, we have $ \;\boldsymbol d_3\cdot\boldsymbol m,_\alpha = \boldsymbol n\cdot\boldsymbol m,_\alpha  =0\; $ and the transverse shear deformations vanishes, since
\begin{equation}\label{f114}
\boldsymbol n_0   \boldsymbol E^e = \boldsymbol n_0 
\big(\boldsymbol Q_e^T\mathrm{Grad}_s\boldsymbol m - \boldsymbol a\big)
 = \big(\boldsymbol n_0 
 \boldsymbol Q_e^T\big)\mathrm{Grad}_s\boldsymbol m
 = \boldsymbol d_3 \big( \boldsymbol m,_\alpha \otimes \boldsymbol a^\alpha\big)
 = (\boldsymbol d_3 \cdot \boldsymbol m,_\alpha ) \boldsymbol a^\alpha
=\boldsymbol 0.
\end{equation}
This shows that the strain shell tensor is a planar tensor in this case, i.e.
\[  
\boldsymbol E^e = E^e_{\alpha\beta} \boldsymbol a^\alpha\otimes \boldsymbol a^\beta \qquad \mathrm{and}\qquad 
\boldsymbol a\boldsymbol E^e = \boldsymbol E^e .
\]
In view of \eqref{f108}, \eqref{f114} and $ \,\boldsymbol b\, \boldsymbol b^* =K\boldsymbol a\, $, we can put the strain-energy density \eqref{f107} in the following reduced form
\begin{equation}\label{f115}
\widetilde W_{\mathrm{shell}} =
\Big(h+ K\, \dfrac{h^3}{12}\Big)   W_{\mathrm{mixt}}\big(  \boldsymbol E^e\big) 
+\,\dfrac{h^3}{12}\, W_{\mathrm{mixt}}\big( \boldsymbol E^e\boldsymbol b 
+ \boldsymbol c \boldsymbol K^e \big)
-2\;\dfrac{h^3}{12}\, W_{\mathrm{mixt}}\big( \boldsymbol E^e\,,\, ( \boldsymbol E^e\boldsymbol b 
+ \boldsymbol c \boldsymbol K^e ) \boldsymbol b^*  \big).
\end{equation}
We see that the right-hand side of \eqref{f115} is a quadratic form of the planar tensors $ \boldsymbol E^e $ and $ \boldsymbol E^e\boldsymbol b + \boldsymbol c\boldsymbol K^e $. Let us express these two tensors in terms of the Koiter shell strain measures $ \boldsymbol \epsilon $ and $ \boldsymbol \rho $.

From \eqref{f112}$ _1 $ and \eqref{f31}$ _1 $ it follows
\begin{equation}\label{f115,5}
\begin{array}{rl}
\boldsymbol \epsilon =&
\dfrac12\,\big[\, (\boldsymbol Q_e^T\mathrm{Grad}_s\boldsymbol m)^T  (\boldsymbol Q_e^T\mathrm{Grad}_s\boldsymbol m) - \boldsymbol a\,\big]
\;=\;
\dfrac12\,\big[\, (\boldsymbol E^e+ \boldsymbol a)^T  (\boldsymbol E^e+ \boldsymbol a) - \boldsymbol a\,\big]
\vspace{6pt}\\
=& \dfrac12\,\big(\boldsymbol E^{e,T}\boldsymbol E^e + \boldsymbol a\boldsymbol E^e+ \boldsymbol E^{e,T}\boldsymbol a) 
\;=\; \dfrac12\,\boldsymbol E^{e,T}\boldsymbol E^e + \mathrm{sym}\big(\boldsymbol a\boldsymbol E^e\big),
\end{array}
\end{equation}
which means
\begin{equation}\label{f116}
\mathrm{sym}\,\boldsymbol E^e\;=\; \boldsymbol \epsilon \,-\,  \dfrac12\,\boldsymbol E^{e,T}\boldsymbol E^e \,.
\end{equation}
Similarly, using \eqref{f112}$ _2\, $, \eqref{f113} and the relation $ \;\boldsymbol Q_e^T \mathrm{Grad}_s\boldsymbol d_3 = \boldsymbol c \boldsymbol K^e -\boldsymbol b \;$ (see \cite[f. (70)]{Birsan-Neff-MMS-2019}), we find
\begin{equation}\label{f117}
\begin{array}{rl}
\boldsymbol \rho =&
- (\boldsymbol Q_e^T\mathrm{Grad}_s\boldsymbol m)^T  (\boldsymbol Q_e^T\mathrm{Grad}_s\boldsymbol d_3) - \boldsymbol b 
\;=\;
- (\boldsymbol E^e+ \boldsymbol a)^T  (\boldsymbol c \boldsymbol K^e -\boldsymbol b) - \boldsymbol b 
\vspace{6pt}\\
=& -\boldsymbol E^{e,T}\boldsymbol c \boldsymbol K^e - \boldsymbol c \boldsymbol K^e + \boldsymbol E^{e,T}\boldsymbol b 
\;=\; -\boldsymbol E^{e,T}\boldsymbol c \boldsymbol K^e - \big( \boldsymbol E^e\boldsymbol b 
+\boldsymbol c \boldsymbol K^e\big) + 2 \big(\mathrm{sym}\,\boldsymbol E^e\big)\boldsymbol b \,.
\end{array}
\end{equation}
Substituting \eqref{f116} in \eqref{f117}, we derive
\begin{equation}\label{f118}
\boldsymbol E^e\boldsymbol b + \boldsymbol c\boldsymbol K^e =
2\, \boldsymbol \epsilon\,\boldsymbol b - \boldsymbol \rho - \boldsymbol E^{e,T}(\boldsymbol E^e\boldsymbol b + \boldsymbol c\boldsymbol K^e).
\end{equation}
With the help of \eqref{f116} and \eqref{f118} we can write now the strain-energy \eqref{f115}  as a function of the strain measures $ \boldsymbol \epsilon $ and $ \boldsymbol \rho $\,: for the first term in 
\eqref{f115} we obtain (from \eqref{f110,5} and \eqref{f115,5})
\begin{equation}\label{f119}
\begin{array}{rl}
W_{\mathrm{Koit}}(\boldsymbol \epsilon) & = \;  \mu\,\|\, \boldsymbol \epsilon\,\|^2  +\,\dfrac{\lambda\,\mu}{\lambda+2\mu}\,\big( \mathrm{tr} \,\boldsymbol \epsilon\big)^2 
\vspace{6pt}\\
& = \;
\mu\,\|\, \mathrm{sym}\,\boldsymbol E^e +   \dfrac12\,\boldsymbol E^{e,T}\boldsymbol E^e \,\|^2  +\,\dfrac{\lambda\,\mu}{\lambda+2\mu}\,\big[ \,\mathrm{tr}\, \big(\mathrm{sym}\,\boldsymbol E^e +   \dfrac12\,\boldsymbol E^{e,T}\boldsymbol E^e\big)\big]^2 .
\end{array}
\end{equation}
Since our model is physically linear (the strain-energy is quadratic in the strain measures) we can neglect the terms in  \eqref{f119} which are more than quadratic in $ \boldsymbol E^e $ and find
\[  
W_{\mathrm{Koit}}(\boldsymbol \epsilon) = \mu\,\|\, \mathrm{sym}\,\boldsymbol E^e \,\|^2  +\,\dfrac{\lambda\,\mu}{\lambda+2\mu}\,\big(\mathrm{tr}\, \boldsymbol E^e \big)^2
\]
i.e.
\begin{equation}\label{f120}
 h W_{\mathrm{Koit}}(\boldsymbol \epsilon) \,= \,h W_{\mathrm{mixt}}(\boldsymbol E^e).
\end{equation}
Thus, the extensional part of our strain-energy density \eqref{f115} coincides in this case with the extensional part of the Koiter model \eqref{f111}.

Similarly, we compute the other two terms of the energy \eqref{f115} and discard the terms which are over-quadratic in the strain measures 
$ \boldsymbol E^e $, $ \boldsymbol K^e $\,: in view of \eqref{f110} and \eqref{f118} we have
\[  
\begin{array}{rl}
W_{\mathrm{Koit}}(\boldsymbol \rho) & =\; W_{\mathrm{mixt}}(\boldsymbol \rho) 
\;=\; 
W_{\mathrm{mixt}}\big(2\, \boldsymbol \epsilon\,\boldsymbol b - (\boldsymbol E^e\boldsymbol b + \boldsymbol c\boldsymbol K^e) - \boldsymbol E^{e,T}(\boldsymbol E^e\boldsymbol b + \boldsymbol c\boldsymbol K^e)\big)
\vspace{6pt}\\
& = \;
W_{\mathrm{mixt}}\big(2\, \boldsymbol \epsilon\,\boldsymbol b - (\boldsymbol E^e\boldsymbol b + \boldsymbol c\boldsymbol K^e)\big)
\vspace{6pt}\\
& = \;
W_{\mathrm{mixt}}\big(\boldsymbol E^e\boldsymbol b + \boldsymbol c\boldsymbol K^e\big)
+ 4\,
W_{\mathrm{mixt}}\big(  \boldsymbol \epsilon\,\boldsymbol b\big)
- 4\,
W_{\mathrm{mixt}}\big( \boldsymbol \epsilon\,\boldsymbol b \,,\, \boldsymbol E^e\boldsymbol b + \boldsymbol c\boldsymbol K^e\big).
\end{array}
\]
It follows
\[  
W_{\mathrm{mixt}}(\boldsymbol E^e\boldsymbol b + \boldsymbol c\boldsymbol K^e) \;=\;  W_{\mathrm{Koit}}(\boldsymbol \rho) -  4\,
W_{\mathrm{mixt}}\big(  \boldsymbol \epsilon\,\boldsymbol b\big)
+ 4\,
W_{\mathrm{mixt}}\big( \boldsymbol \epsilon\,\boldsymbol b \,,\, \boldsymbol E^e\boldsymbol b + \boldsymbol c\boldsymbol K^e\big)
\]
and inserting \eqref{f118} here we find for the second term in the energy \eqref{f115}:
\begin{equation}\label{f121}
\begin{array}{rl}
W_{\mathrm{mixt}}(\boldsymbol E^e\boldsymbol b + \boldsymbol c\boldsymbol K^e) & =  \; W_{\mathrm{Koit}}(\boldsymbol \rho) -  4\,
W_{\mathrm{mixt}}\big(  \boldsymbol \epsilon\,\boldsymbol b\big)
+ 4\,
W_{\mathrm{mixt}}\big( \boldsymbol \epsilon\,\boldsymbol b \,,\, 2\, \boldsymbol \epsilon\,\boldsymbol b - \boldsymbol \rho \big)
\vspace{6pt}\\
& = \;
 W_{\mathrm{Koit}}(\boldsymbol \rho) +  4\,
 W_{\mathrm{mixt}}\big(  \boldsymbol \epsilon\,\boldsymbol b\big)
 - 4\,
 W_{\mathrm{mixt}}\big( \boldsymbol \epsilon\,\boldsymbol b \,,\,  \boldsymbol \rho \big).
\end{array}
\end{equation}
For the last term in \eqref{f115} we write with the help of \eqref{f118}:
\begin{equation}\label{f122}
(\boldsymbol E^e\boldsymbol b + \boldsymbol c\boldsymbol K^e)\boldsymbol b^* \;=\; 2K\, \boldsymbol \epsilon - \boldsymbol \rho\, \boldsymbol b^*  - 
\boldsymbol E^{e,T}(\boldsymbol E^e\boldsymbol b + \boldsymbol c\boldsymbol K^e)
\boldsymbol b^* 
\end{equation}
and derive from \eqref{f116} and \eqref{f122}
\begin{equation}\label{f123}
\begin{array}{l}
W_{\mathrm{mixt}}\big(\boldsymbol E^e, (\boldsymbol E^e\boldsymbol b + \boldsymbol c\boldsymbol K^e)\boldsymbol b^*\big) \; =  \;
W_{\mathrm{mixt}}\big(\mathrm{sym}\,\boldsymbol E^e,\,  2K\, \boldsymbol \epsilon - \boldsymbol \rho\, \boldsymbol b^*\big)
\vspace{6pt}\\
\qquad\qquad  = \;
 W_{\mathrm{mixt}}\big(\, \boldsymbol\epsilon\,,\,  2K\, \boldsymbol \epsilon - \boldsymbol \rho\, \boldsymbol b^*\big) 
 \;=\; 2K\, W_{\mathrm{Koit}}( \boldsymbol\epsilon ) -  W_{\mathrm{mixt}}\big(\, \boldsymbol\epsilon\,,\, \boldsymbol \rho\, \boldsymbol b^*\big),
\end{array}
\end{equation}
We substitute \eqref{f120}, \eqref{f121} and \eqref{f123} into \eqref{f115} and obtain
\[  
\begin{array}{rl}
\widetilde W_{\mathrm{shell}}(\boldsymbol \epsilon,\boldsymbol \rho) \;  = &
\Big(h+ K\, \dfrac{h^3}{12}\Big)   W_{\mathrm{Koit}}\big(  \boldsymbol \epsilon\big) 
+\,\dfrac{h^3}{12}\,\Big( W_{\mathrm{Koit}}( \boldsymbol \rho )
+ 4 \, W_{\mathrm{mixt}}\big( \boldsymbol \epsilon\, \boldsymbol b\big)
-  4 \, W_{\mathrm{mixt}}\big( \boldsymbol \epsilon\, \boldsymbol b
\,,\, \boldsymbol \rho\big)
\Big)
\vspace{6pt}\\
&
-2\;\dfrac{h^3}{12} \,\Big( 2K\, W_{\mathrm{Koit}}( \boldsymbol \epsilon )
- W_{\mathrm{mixt}}\big( \boldsymbol \epsilon
\,,\, \boldsymbol \rho\, \boldsymbol b^*\big)
\Big),
\end{array}
\]
which can be written in view of \eqref{f110} in the form
\begin{equation}\label{f124}
\begin{array}{c}
\widetilde W_{\mathrm{shell}}(\boldsymbol \epsilon,\boldsymbol \rho) \,  = \,
h\,   W_{\mathrm{Koit}}\big(  \boldsymbol \epsilon\big) 
+\,\dfrac{h^3}{12}\, W_{\mathrm{Koit}}( \boldsymbol \rho )
+\,\dfrac{h^3}{12}\, \Big[
 4 \, W_{\mathrm{mixt}}\big( \boldsymbol \epsilon\, \boldsymbol b
 \,,\, \boldsymbol \epsilon\, \boldsymbol b - \boldsymbol \rho
 \big)
- \, W_{\mathrm{mixt}}\big( \boldsymbol \epsilon\,,\, 3K\,\boldsymbol \epsilon - 2 \boldsymbol\rho \,\boldsymbol b^*\big)
\Big].
\end{array}
\end{equation}
The terms in the square brackets in \eqref{f124} involve the initial curvature of the shell through the tensor $ \,\boldsymbol b \,$, the cofactor $ \,\boldsymbol b^* = 2H  \boldsymbol a - \boldsymbol b\, $ and the determinant $\, K=\mathrm{det}\, \boldsymbol b\, $ (Gau\ss{} curvature). These terms vanish  in the case of plates (since $ \boldsymbol b = \boldsymbol 0 $); moreover,  they
can be neglected also for sufficiently thin shells, provided the midsurface strain is small.
We note that the corresponding terms in the classical shell theory have been neglected using
similar arguments, see the discussion about the term $ W_3 $ in \cite[f. (57)]{Steigmann13}. Finally, if we retain only the leading extensional and bending terms in \eqref{f124}, we obtained the reduced classical form
\begin{equation}\label{f125}
\widetilde W_{\mathrm{shell}}(\boldsymbol \epsilon,\boldsymbol \rho) 
\,  = \,
h\,   W_{\mathrm{Koit}}\big(  \boldsymbol \epsilon\big) 
+\,\dfrac{h^3}{12}\, W_{\mathrm{Koit}}( \boldsymbol \rho )\,,
\end{equation}
in accordance with the Koiter energy density \eqref{f111}.

In conclusion, our model can be regarded as a generalization of the classical Koiter model in the framework of 6-parameter shell theory.

\bigskip\bigskip
\noindent
\small{\textbf{Acknowledgements}\quad
	This research has been funded by the Deutsche Forschungsgemeinschaft (DFG, German Research Foundation) -- Project no. 415894848.

\bibliographystyle{plain} 

\bibliography{literatur_Birsan}

\end{document}